\documentclass[a4paper,12pt]{amsart}
\usepackage[foot]{amsaddr}
\usepackage[margin=2.7cm]{geometry}
\usepackage{systeme}
\usepackage{graphicx}
\usepackage[]{units}
\usepackage{color}
\usepackage{mathrsfs}
\usepackage{bm}
\usepackage{float}
\usepackage{amssymb}
\usepackage{amsthm}
\usepackage{amsmath}
\usepackage{resizegather}
\usepackage{easyReview}
\usepackage{hyperref} 
\usepackage{color}
\usepackage{ulem}

\newtheorem{theorem}{Theorem}[section]

\theoremstyle{definition}

\theoremstyle{remark}
\newtheorem{remark}[theorem]{Remark}

\numberwithin{equation}{section}

\definecolor{nverde}{RGB}{0,61,0}

\title[POD-Galerkin ROM combined Navier-Stokes transport equations]{POD-Galerkin reduced order methods for combined Navier-Stokes transport equations based on a hybrid FV-FE solver}

\date{}

\begin{document}
	\author{S. Busto\textsuperscript{1,2,*}}
	\address{\textsuperscript{1} Departamento de Matem\'atica Aplicada, Universidade de Santiago de Compostela. ES-15782 Santiago de Compostela, Spain.}
	\address{\textsuperscript{2} Department of Mathematics, Unit\`a INdAM, University of Trento. IT-38100 Trento, Italy}
	\thanks{\textsuperscript{*}Corresponding Author.}
	\email{saray.busto@usc.es}
	
	\author{G. Stabile\textsuperscript{3}}
	\address{\textsuperscript{3} SISSA, International School for Advanced Studies, Mathematics Area, mathLab, Trieste, 34136, Italy.}
	\email{gstabile@sissa.it}
	
	\author{G. Rozza\textsuperscript{3}}
	\email{grozza@sissa.it}
	
	\author{M.E. V\'azquez-Cend\'on\textsuperscript{1}}
	\email{elena.vazquez.cendon@usc.es}
	
	\keywords{proper orthogonal decomposition, projection hibrid finite volume-finite element method, Poisson pressure equation, Navier-Stokes equations, transport equation}

\begin{abstract}
The purpose of this work is to introduce a novel POD-Galerkin strategy for the hybrid finite volume/finite element solver introduced in \cite{BFSV14} and \cite{BFTVC17}. The interest is into the incompressible Navier-Stokes equations coupled with an additional transport equation.
The full order model employed in this article makes use of staggered meshes. This feature will be conveyed to the reduced order model leading to the definition of reduced basis spaces in both meshes.
 The reduced order model presented herein accounts for velocity, pressure, and a transport-related variable. The pressure term at both the full order and the reduced order level is reconstructed making use of a projection method. More precisely, a Poisson equation for pressure is considered within the reduced order model. Results are verified against three-dimensional manufactured test cases. Moreover a modified version of the classical cavity test benchmark including the transport of a species is analysed.
\end{abstract}

	\maketitle

% % % % % % % % % % % % % % % % % % % % % % % % % % % % % %
%                    Introduction                         %
% % % % % % % % % % % % % % % % % % % % % % % % % % % % % %
\section{Introduction}\label{sec:intro}

During the last decades, numerical methods to solve partial differential equations (PDEs) have became more efficient and reliable and it is now possible to find a large variety of different solvers using diverse discretization methods. Numerical methods transform the original set of PDEs into a possibly very large system of algebraic equations and there are still many situations where this operation, using standard discretization techniques (Finite Difference Methods (FDM), Finite Element Methods (FEM), Finite Volume Methods (FVM),...), becomes not feasible. Such situations occur any time that we need to solve the problem in a parametrized setting in order to perform optimization or uncertainty quantification or when we need a reduced computational cost  as in real-time control. A possible way to deal with the former situations, at a reasonable computational cost, is given by the reduced order methodology (ROM). In ROM the high dimensional full order model is replaced with a surrogate reduced model that has a much smaller dimension and therefore is much cheaper to solve (\cite{Pet89,IR98,NE94,HRS16,quarteroniRB2016,ChinestaEnc2017,bennerParSys}). In this manuscript we focus our attention on projection-based reduced order models and in particular on POD-Galerkin methods (\cite{IL98a,DKKO91}). They have been successfully applied in a wide variety of different settings ranging from computational fluid dynamics to structural mechanics. The scope of this paper is the development of a POD-Galerkin reduced order method for the unsteady incompressible Navier-Stokes equations coupled with a transport equation. The full order solver, on which this work is based, is the hybrid FV/FE solver presented in \cite{BFSV14} and \cite{BFTVC17}. It exploits the advantages given by both full order discretization techniques. 

The main novelties of the present work consist into the development of a reduced order for unsteady flows starting from an hybrid FV-FE solver on a dual mesh structure and into its coupling with an additional transport equation. The reduced basis for the generation of the POD spaces are in fact defined on two different meshes and, in order to map the basis functions between the finite element mesh and the finite volume one, we exploit the structure of the full order solver. Standard POD-Galerkin methods developed for FE and FV discretization have been re-adapted to be used in the hybrid finite volume-finite element framework introduced in \cite{BFSV14,BFTVC17,BBFVC18}.

The paper is organized as follows. In \autoref{sec:FOM} the hybrid FV/FE solver is explained and all the main important features, which are relevant in a reduced order modeling setting, are recalled. In \autoref{sec:ROM} the reduced order model is introduced with all the relevant details and remarks. In \autoref{sec:numerical_results} two different numerical tests are presented. The tests consist into a manufactured three-dimensional fluid dynamic problem and into a modified three-dimensional lid-driven cavity problem combined with a species transport equation. Finally, in \autoref{sec:conclusions} some conclusions and perspectives for future works are provided.

\vspace{0.2cm}

% % % % % % % % % % % % % % % % % % % % % % % % % % % % % %
%                         FOM                             %
% % % % % % % % % % % % % % % % % % % % % % % % % % % % % %
\section{The Full Order Model}\label{sec:FOM}

The model for incompressible newtonian fluids is enlarged with a transport equation. Hence, the system of equations to be solved written in conservative form results 
\begin{align}
&\mathrm{div}\mathbf{w}_\mathbf{u}  =0, \label{eq:masa}\\
&\frac{\partial \mathbf{w}_\mathbf{u} }{\partial t} + { \textrm{div} \mathbf{\mathcal{
			F}^{w_u}}(\mathbf{w}_\mathbf{u} )} +{ \mathrm{grad}\, \pi }
-\mathrm{div} \, \tau=\mathbf{f}_{\mathbf{u}}, \label{eq:momentos} \\
% Species
&\frac{\partial w_{\mathrm{y}}}{\partial t}+\mathrm{div} \mathcal{F}^{w_{\mathrm{y}}}\left(w_{\mathrm{y}},\mathbf{u}\right) 
-\mathrm{div}\left[\rho \mathcal{D}\mathrm{grad}\left( \frac{1}{\rho}w_{\mathrm{y}}\right) \right]=0,\label{eq:especies2}
\end{align}
where $\rho$ is the density, $ \mathbf{w_u}:=\rho\mathbf{u}$ is the linear momentum density, $w_{\mathrm{y}}$ is the unknown related to the transport equation that can be, for example, the conservative variable  related to a species $\mathrm{y}$,
$\pi$ is the pressure perturbation, $\tau$ is the viscous part of Cauchy stress tensor, $\mathcal{D}$ is the viscous coefficient related to the transport equation, $\mathbf{f}_{\mathbf{u}}$ is an arbitrary source term for the momentum equation and
\begin{equation}
\mathcal{F}(\mathbf{w})= \left( \mathcal{F}^{\mathbf{w_u}}\left( \mathbf{w}_\mathbf{u} \right) , \mathcal{F}^{w_{\mathrm{y}}}\left(w_{\mathrm{y}},\mathbf{u}\right) \right)^{T},
\end{equation}
denotes the complete flux tensor whose three components read
\begin{equation}
\mathcal{F}=\left( \mathcal{F}_1 | \mathcal{F}_2 | \mathcal{F}_3 \right)_{\left( 3+1\right) \times 3}, \quad \mathcal{F}_i\left(\mathbf{w}\right) =
u_i \mathbf{w},\quad i=1,2,3.
\end{equation}
The numerical discretization of the above system, \eqref{eq:masa}-\eqref{eq:especies2}, is performed by considering the projection hybrid finite volume-finite element method presented in \cite{BFSV14}, \cite{BFTVC17} and \cite{BBFVC18}. In the following section we summarize the main features of the aforementioned methodology.
\subsection{Numerical discretization}\hfill\\

\begin{figure}
	\centering
	\includegraphics[width=0.45\linewidth]{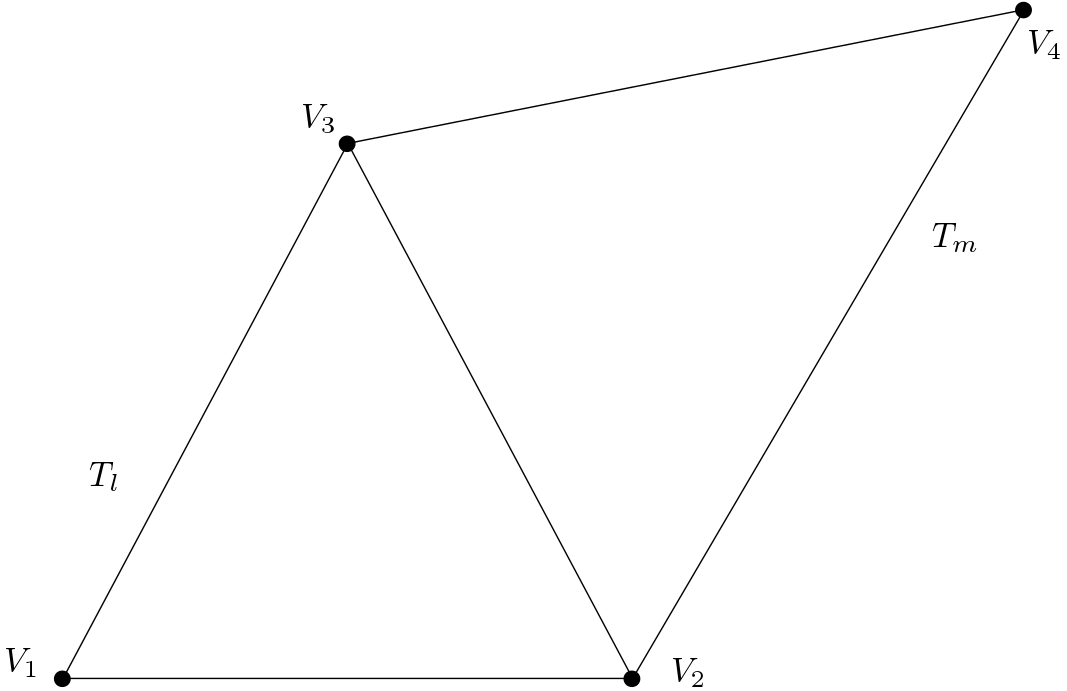}
	\includegraphics[width=0.45\linewidth]{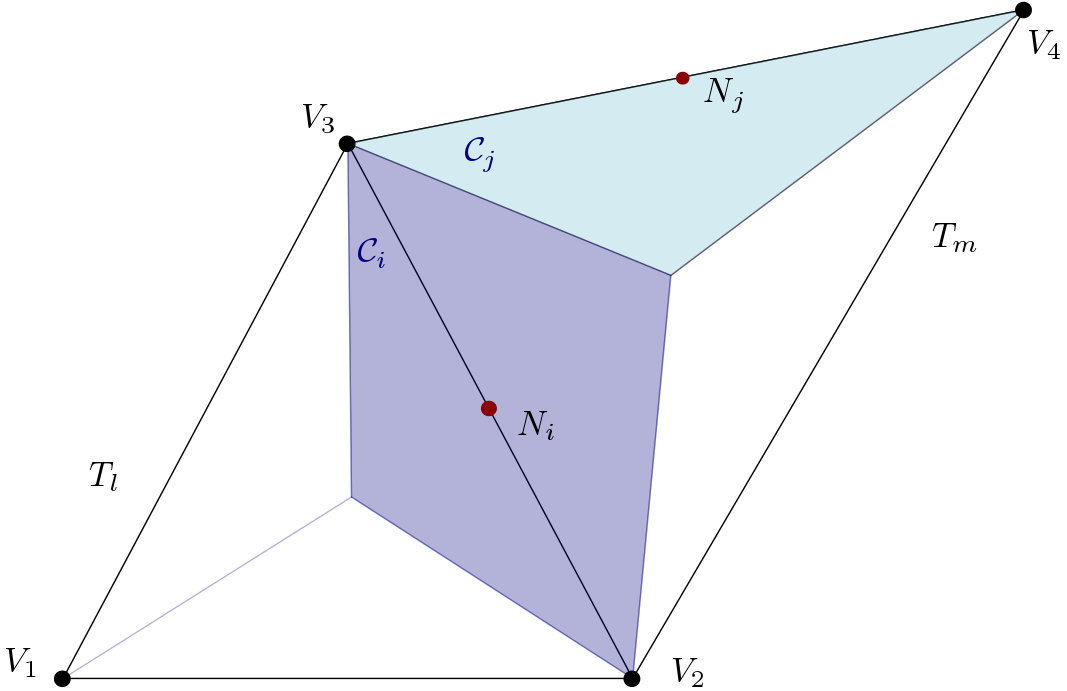}
	\caption{Construction of a dual mesh in 2D. Left: Finite elements
		of the primal mesh (black). Right: interior finite volume $C_{i}$ (purple) of the dual mesh related to elements $T_{l}$ and $T_{m}$ of the primal mesh and boundary finite volume $C_{j}$ (blue) related to the finite element $T_{m}$.}
	\label{fig:descentradosnuevoef}
\end{figure}

A two-stage in time discretization scheme is considered:
in order to get the solution
at time $t^{n+1}$, we use the previously obtained approximations
$\mathbf{W}^n$ of the conservative variables $\mathbf{w}(x,y,z,t^n)$,
$\mathbf{U}^n$ of the velocity $\mathbf{u}(x,y,z,t^n)$
and $\pi^n$ of the pressure perturbation $\pi(x,y,z,t^n)$, and
compute $\mathbf{W}^{n+1}$ and $\pi^{n+1}$ from the following system
of equations:
\begin{eqnarray}
\frac{1 }{\Delta t}\left( \widetilde{\mathbf{W}}^{n+1}_{\mathbf{u}}-\mathbf{W}^{n}_{\mathbf{u}}\right)   +  \mathrm{div} \mathbf{\mathcal{F}}^{\mathbf{w}_{\mathbf{u}}}(
\mathbf{W}^{n}_{\mathbf{u}})
+  \nabla \pi^{n} - \mathrm{div}( \tau^n)=\mathbf{f}_{\mathbf{u}}^{n} , \label{eq:tidenincrec} \\
\frac{ 1 }{\Delta t} \left( \mathbf{W}^{n+1}_{\mathbf{u}}-  \mathbf{\widetilde{W}}^{n+1}_{\mathbf{u}}\right)  + 
\nabla (  \pi^{n+1} - \pi^{n}) =0, \label{eq:tideincre2} \\
\mathrm{div}\mathbf{W}^{n+1}_{\mathbf{u}}=0, \label{eq:qincrec}\\
% Especies
\;\;\; \frac{1}{\Delta t}\!\left( W_{\mathrm{y}}^{n+1}- W_{\mathrm{y}}^n\right) +\mathrm{div}\mathcal{F}^{w_{\mathrm{y}}} \left(W_{\mathrm{y}}^n,\mathbf{U}^{n} \right) -\mathrm{div}\left[ \rho\mathcal{D} \,\mathrm{grad}\, \mathrm{y}^n \right] = 0.\label{eq:especies_discret}
\end{eqnarray}\\
Concerning the discretization of mass conservation and momentum equations,
by adding equations \eqref{eq:tidenincrec}-\eqref{eq:tideincre2}, we easily see that the scheme is
actually implicit for the pressure term. However, the equations above
show that the pressure and the velocity can be solved in three uncoupled stages:
\begin{itemize}
	\item Transport-diffusion stage. Equations \eqref{eq:tidenincrec} and \eqref{eq:especies_discret} are explicitly solved by considering a finite volume scheme. We notice that, in general, the intermediate approximation of the linear momentum computed, $\widetilde{\mathbf{W}}^{n+1}_{\mathbf{u}}$, does not  satisfy the divergence free condition	\eqref{eq:qincrec}.
	\item Projection stage. It is a implicit stage in which the coupled equations \eqref{eq:tideincre2} and \eqref{eq:qincrec} are solved with a finite element method to obtain the pressure correction $\delta^{n+1}:=\pi^{n+1}-\pi^{n}$.
	\item Post-projection stage. The intermediate approximation for the
	linear momentum is updated with the pressure correction providing
	the final approximation $\mathbf{W}^{n+1}_{\mathbf{u}}$. Finally, $\pi^{n+1}$ is recovered.
\end{itemize}
For the spatial discretization of the domain we consider an unstructured tetrahedral finite element mesh from which a dual finite volume mesh of the face type is build.  The pressure is approximated at the vertex of the original tetrahedral mesh whereas the
conservative variables are computed in the nodes of the dual mesh. 
Figure \ref{fig:descentradosnuevoef} depicts the construction of the dual mesh in 2D, further details on both the 2D and the 3D cases can be found in \cite{BDDV98} and \cite{BFSV14}.
Denoting by $C_{i}$ a cell of the dual mesh, $N_{i}$ its node, $\Gamma_{i}$ the boundary, $\boldsymbol{\tilde{\eta}}_{i}$ the outward unit normal at the boundary and $\left|C_i\right|$ the volume of the cell, the discretization of the transport-diffusion equations, \eqref{eq:tidenincrec} and \eqref{eq:especies_discret}, results
\begin{gather}
\displaystyle \frac{1}{\Delta t}\left( \widetilde{\mathbf{W}}_{\mathbf{u},\, i}^{n+1} - \mathbf{W}_{\mathbf{u},\ i}^{n}\right)  + \frac{1}{\left| C_{i} \right|} \displaystyle \int_{\Gamma_i}\mathcal{F}^{\mathbf{w}_{\mathbf{u}}}\left({\mathbf{W}}^n_{\mathbf{u}}\right){\boldsymbol{\tilde{\eta}}_{i}}\,\mathrm{dS} 
+ \frac{1}{\left| C_{i} \right|} \displaystyle \int_{C_i}\nabla\pi^n \, \mathrm{dV}  -\frac{1}{\left| C_{i} \right|} \displaystyle\int_{\Gamma_i} \tau^{n}\boldsymbol{\tilde{\eta}}_{i} \,\mathrm{dS}  =
\frac{1}{\left|C_{i}\right|} \int_{C_i} \mathbf{f}_{\mathbf{u}} \, \mathrm{dV}, \label{eq:wu_discret}\\[0.4cm]
\displaystyle \frac{1}{\Delta t}\left( W_{\mathrm{y},\, i}^{n+1} - W_{\mathrm{y},\, i}^n\right)  + \frac{1}{\left| C_{i} \right|} \displaystyle \int_{\Gamma_i} \mathcal{F}^{w_{\mathrm{y}}}\left(W_{\mathrm{y}}^n,\mathbf{U}^{n} \right)\boldsymbol{\tilde{\eta}}_{i} \, \mathrm{dS}  
- \frac{1}{\left| C_{i} \right|} \displaystyle\int_{\Gamma_i}   \mathcal{D}\, \mathrm{grad}\,   W_{\mathrm{y}}^n \, \boldsymbol{\tilde{\eta}}_{i}\,\mathrm{dS} = 0,\label{eq:y_discret}
\end{gather}
The flux term on the above expressions is computed by considering the Rusanov scheme. Moreover, to obtain a second order in space scheme the CVC Kolgan-type method is used (see \cite{CVC10}, \cite{CV12} and \cite{BFTVC17}). Second order in both space and time is achieved by applying LADER methodology (see \cite{BFTVC17} and \cite{BBFVC18}). This last technique also affects the computation of the viscous term in which a Galerkin approach is consider to approximate the spacial derivatives. Finally, the pressure term is obtained by averaging its values at the three vertex of each face and the barycenter of the tetrahedra to which the face belongs.
Within the projection stage, the pressure is computed using a standard finite element method. The incremental projection method presented in \cite{Guer06} is adapted to solve \eqref{eq:tideincre2}-\eqref{eq:qincrec}
obtaining the following weak problem:
{\textit{Find $\delta^{n+1}\in V_0:=\left\{ z\in H^1(\Omega): \int_{\Omega}z=0\right\}$ verifying }}
\begin{equation}
\displaystyle\int_{\Omega} \nabla \delta^{n+1}\cdot \nabla z  \, \mathrm{dV}= \displaystyle \frac{1}{\Delta t}
\int_{\Omega}\mathbf{\widetilde{W}}^{n+1} \cdot \nabla z  \, \mathrm{dV}  - \displaystyle \frac{1}{\Delta t} \int_{\partial \Omega} G^{n+1}z  \, \mathrm{dS}\quad \forall z\in V_0, 
\end{equation}
where $\delta^{n+1}:= \pi^{n+1}-\pi^{n}$ and $G^{n+1}= \mathbf{\widetilde{W}}^{n+1}_{\mathbf{u}} \cdot \boldsymbol{\eta}$ (see \cite{BFSV14} for further details).

Finally, at the post-projection stage, $\mathbf{W}_{\mathbf{u}}^{n+1}$ is updated by using $\delta_{i}^{n+1}$:
\begin{equation}
\mathbf{W}_{\mathbf{u}}^{n+1}=\widetilde{\mathbf{W}}_{\mathbf{u}}^{n+1}+\Delta t\, \mathrm{grad}\, \delta_{i}^{n+1}.
\end{equation}

\vspace{0.2cm}
% % % % % % % % % % % % % % % % % % % % % % % % % % % % % %
%                         ROM                             %
% % % % % % % % % % % % % % % % % % % % % % % % % % % % % %
\section{The Reduced Order Model}\label{sec:ROM}

In this section, projection based ROM techniques are adapted to be used in the hybrid finite element-finite volume framework already presented for the FOM. Projection based ROMs haven been used to deal with different type of applications starting from both a finite element discretization \cite{Ball15,quarteroniRB2016,HRS16} and a finite volume one \cite{DHO2012,HOR08,SR17}.

To this end, two different POD spaces are considered. The first one is related to the linear momentum and provides the solution at the nodes of the finite volume mesh. On the other hand, the POD basis associated to the pressure are computed on the vertex of the primal mesh. Moreover, when considering also the transport equation a third POD space for the transport unknown are generated.

\subsection{The POD spaces}\label{sec:pod_spaces} \hfill\\

We start by introducing the computation of the POD space related to the linear momentum. The main assumption in the POD-ROM technique is that the approximated solution $\mathbf{w}_\mathbf{u}^{n}$ can be expressed as linear combination of spatial modes $\boldsymbol{\varphi}_i(\mathbf{x})$ multiplied by scalar temporal coefficients $a_i(t)$ (\cite{NE94,quarteroniRB2016,HRS16}), that is,
\begin{equation}
\mathbf{w}_\mathbf{u}(\mathbf{x},t^{n})=\sum_{i=1}^{N} a_i(t^{n}) \boldsymbol{\varphi}_i(\mathbf{x}),\label{eq:lc_velocity}
\end{equation}
where $N$ is the cardinality of the basis of the spatial modes.
To compute the basis elements we will use POD. This methodology allows us to select the most energetic modes so that we maximize the information contained in the basis.

Let us consider a training set 
$\mathcal{K}=\left\lbrace t^{1},\dots,t^{N_{s}} \right\rbrace\subset [0,T]$
of discrete time instants taken in the simulation time interval and its corresponding set of snapshots
$\left\lbrace \mathbf{w}_{\mathbf{u}}^{n}\left(\mathbf{x},t^{n}\right), \, t^{n}\in \mathcal{K} \right\rbrace$. 
Then, the elements of the POD basis are chosen to minimize the difference between the snapshots and their projection in $L^2$ norm
\begin{equation}
X^{\mathrm{POD}}_{N_r} =  \textrm{arg min} \frac{1}{N_s} \sum_{n=1}^{N_s}\left\| \mathbf{w}_{\mathbf{u}}^{n}\left(\mathbf{x}\right) - \sum_{i=1}^{N}\left\langle \mathbf{w}_{\mathbf{u}}^{n}\left(\mathbf{x}\right),\boldsymbol{\varphi}_i(\mathbf{x}) \right\rangle_{L^2}  \boldsymbol{\varphi}_i(\mathbf{x})   \right\|^{2}_{\mathrm{L}^2}
\label{eq:minimizationpb_u}\end{equation}
with $\left\langle \boldsymbol{\varphi}_{i}(\mathbf{x}),\boldsymbol{\varphi}_{j}(\mathbf{x}) \right\rangle_{L^2\left(\Omega \right)} =\delta_{ij}$ for all $i,j\in\left\lbrace1,\dots,N\right\rbrace$. The minimization problem  can be solved following \cite{KV02}
where it has been shown that \eqref{eq:minimizationpb_u} is equivalent to the eigenvalue problem
\begin{equation}
\bm{\mathcal{C}} \bm{\xi}=\bm{\lambda}\bm{\xi},
\end{equation}
with
\begin{equation}
\mathcal{C}_{jk}=\frac{1}{N_s}\left\langle\mathbf{w}_{\mathbf{u}}^{j}\left(\mathbf{x}\right) ,\mathbf{w}_{\mathbf{u}}^{k}\left(\mathbf{x}\right)  \right\rangle_{L^2}
=\frac{1}{N_s} \sum_{i=1}^{N_{nod}} \left[\left|C_{i}\right|  \sum_{l=1}^{3} w_{l}^{j}\left(\mathbf{x}_i\right)w_{l}^{k}\left(\mathbf{x}_i\right) \right]
\end{equation}
being the correlation matrix of the snapshots, $\lambda_{i}$ the eigenvalues, $\bm{\xi}_{i}$ the corresponding eigenvectors and $N_{nod}$ the number of elements of the FV mesh.
Since this methodology sorts the eigenvalues in descending order, the first modes retain most of the energy present in the original solutions. Defining the cumulative energy as
\begin{equation*}E_{\lambda}:=\sum_{i=1}^{N_{s}} \lambda_i,\end{equation*}
the number of elements of the POD basis can be set as the minimum $N$ such that $\sum_{i=1}^{N} \left( \frac{\lambda_i}{E_{\lambda}}\right) $ is lower or equal to a desired fixed bound $\kappa_{\mathbf{w}_{\mathbf{u}}}$. Next, the elements of the basis are computed as
\begin{equation}
\boldsymbol{\varphi}_i(\mathbf{x})=\frac{1}{\sqrt{\lambda_i}}\sum_{n=1}^{N_s} \xi_{in}\mathbf{w}_{\mathbf{u}}^{n}(\mathbf{x}), \quad i=1,\dots, N,
\end{equation}
and normalized. Finally, the POD space is given by
\begin{equation}
X^{\mathrm{POD}}_{N} := \displaystyle{\textrm{span}_{i=1,\dots,N} } \left\lbrace\boldsymbol{\varphi}_i\right\rbrace.
\end{equation}

\begin{remark}[Source term]
	In some cases, for example when we are working with analytical tests, we will be interested in including a source term in the momentum equation. 
	If there exists a linear relation between the source term and the velocity field, we can compute a basis for it considering the eigenvalue problem already solved for the velocity. Then, the elements of the source term basis would be computed as
	\begin{equation}
	\boldsymbol{\varsigma}_i(\mathbf{x})=\frac{1}{\sqrt{\lambda_i} \left\| \mathbf{\varphi}_i (\mathbf{x}) \right\|_{L^2} }\sum_{n=1}^{N_s} \xi_i \mathbf{f}_{\mathbf{u}}^n(\mathbf{x}).
	\end{equation}	
	Therefore, the source term results
	\begin{equation}
	\mathbf{f_u}\left( \mathbf{x},t\right) = \sum_{n=1}^{N} a_i(t)\boldsymbol{\varsigma}_i(\mathbf{x}).\label{eq:lc_source}
	\end{equation}
\end{remark}	

\begin{remark}
	In the finite element framework, since the natural functional space for velocity is $\mathcal{H}^{1}$, to compute its correlation matrix, $\mathcal{C}$, the $H^{1}$ norm is used whereas the $L^{2}$ norm is employed only to compute the correlation matrix for pressure. However, in the finite volume framework the $L^{2}$ norm is also used for the velocities. Since the velocity belongs to a discontinuous space and so, the computation of the gradient to evaluate the $H^{1}$ norm will introduce further discretisation error (see \cite{SR17}). Moreover, the $L^{2}$ norm has a direct physical meaning being directly correlated to the kinetic energy of the system.
\end{remark}

Let us assume that the pressure can be expressed as a linear combination of the elements of a POD basis
\begin{equation}
\pi(\mathbf{x},t^{n})=\sum_{i=1}^{N_{pi}} b_i(t) \psi_i(\mathbf{x}), \quad
X^{\mathrm{POD}_{\pi}}_{N_{\pi}} := \displaystyle{\textrm{span}_{i=1,\dots,N_{\pi}} } \left\lbrace \psi_i\right\rbrace. \label{eq:lc_pressure}
\end{equation}
The  modes of the above basis can be computed similarly to the ones for the linear momentum. Nonetheless, it is important to notice that the $L^{2}$ products involving the pressure are now computed in the finite element framework. More precisely, we use a quadrature rule on the barycentre of the edges of the tetrahedra, so that the correlation matrix is given by
\begin{equation}
\mathcal{C}_{jk}^{\pi}=\frac{1}{N_s}\left\langle\!\!\left\langle\pi^{j}\left(\mathbf{x}\right) ,\pi^{k}\left(\mathbf{x}\right)  \right\rangle\!\!\right\rangle_{L^2}
=\frac{1}{N_s}\sum_{i=1}^{N_{nel}} \left(\left|T_{i}\right|  \sum_{l=1}^{6} \omega_{l} \pi^{j}\left(\mathbf{x}_{il}\right)\pi^{k}\left(\mathbf{x}_{il}\right) \right)
\end{equation}
with $\mathbf{x}_{il}$ the barycentre of the $l-th$ edge of the tetrahedra $T_{i}$, $\omega_{l}=\frac{1}{6}$ the weights and $N_{nel}$ the number of elements of the FE mesh.
Besides, we are assuming that the pressure basis is  independent to the one corresponding to the linear momentum. Accordingly, the dimension of the diverse basis may be different.

Eventually, since the basis functions for the variable related to the transport equation may also be linearly independent from the ones obtained for the linear momentum, a new spectral problem will be defined in the FV framework. As a result, we obtain the POD basis space for $w_\mathrm{y}$:
\begin{equation}
X^{\mathrm{POD}_{\mathrm{y}}}_{N_{\mathrm{y}}} := \displaystyle{\textrm{span}_{i=1,\dots,N_{\mathrm{y}}} } \left\lbrace \chi_i\right\rbrace, \quad w_\mathrm{y}(\mathbf{x},t^{n})=\sum_{i=1}^{N_{\mathrm{y}}} c_i(t) \chi_i(\mathbf{x}). \label{eq:lc_species} 
\end{equation}

% % % % % % % % % % % % % % % % % % % % % % % % % % % % % %
% % % % % % % % % % % % % % % % % % % % % % % % % % % % % %
\subsection{Galerkin projection}\hfill\\

Once the basis functions for the POD spaces are obtained we need to compute the unknown vectors of coefficients $\mathbf{a}$, $\mathbf{b}$ and $\mathbf{c}$. To this end, we perform a Galerkin projection of the governing equations onto the POD reduced basis spaces and we solve the resulting algebraic system.

% % % % % % % % % % % % % % % % % % % % % % % % % % % % % %
\subsubsection{Momentum equation}\hfill\\

Regarding the momentum equation, the resulting dynamical system reads
\begin{equation}
\dot{\mathbf{a}} = -\mathbf{a}^{T} \mathbf{C}\mathbf{a}+ \mathbf{B} \mathbf{a} - \mathbf{K} \mathbf{b}  + \mathbf{F} \mathbf{a},\label{eq:dyn_sys_w}
\end{equation}
where
\begin{gather}
B_{ij}=\langle \boldsymbol{\varphi}_{i}, \frac{\mu}{\rho} \Delta \boldsymbol{\varphi}_{j} \rangle_{L^{2}},\quad
C_{ijk}= \langle \boldsymbol{\varphi}_{i}, \frac{1}{\rho}\mathrm{div}\,\left(  \boldsymbol{\varphi}_{j}\otimes \boldsymbol{\varphi}_{k}\right)  \rangle_{L^{2}},\notag\\
K_{ij}=  \left< \boldsymbol{\varphi}_{i}, \mathrm{grad}\,\psi_{j} \right>_{L^{2}},\quad
F_{ij}=  \left< \boldsymbol{\varphi}_{i},  \boldsymbol{\varsigma}_{j} \right>_{L^{2}}.
\end{gather}

The computation of the gradients of the basis functions involved in the former products are performed considering the methodology already introduced for the computation of the spatial derivatives of the viscous term on the FOM. Thus we take advantage of the dual mesh structure. We will further detail it in the case of the convective term.

Due to the non-linearity of the convective term, the computation of its related matrix is more complex than the ones of the remaining terms and requires for the calculation and storage of a third-order tensor. Substituting $\mathbf{w}_{\mathbf{u}}$ by the corresponding linear combination of the POD modes and projecting the convection term, we get
\begin{gather}
\langle \boldsymbol{\varphi}_{i}, \frac{1}{\rho}  \mathrm{div}\, \left( \mathbf{w}_{\mathbf{u}}\otimes \mathbf{w}_{\mathbf{u}}\right) \rangle_{L^{2}} =
\langle \boldsymbol{\varphi}_{i}, \frac{1}{\rho} \mathrm{div}\left(  \sum_{j=1}^{N} a_j \boldsymbol{\varphi}_j \otimes \sum_{k=1}^{N} a_k \boldsymbol{\varphi}_k \right) \rangle_{L^{2}} 
=\mathbf{a}^{T} \mathbf{C}\mathbf{a}.
\end{gather}
In order to develop an efficient approach complete consistency with the FOM has not been respected within the computation of the above term.
Moreover, the use of the former approach simplifies the calculations reducing the computational cost of the method. Still its derivation is consistent with the physical model. 

Moreover, the tools employed in the proposed approach have already been successfully used to compute the contribution of the advection term on reconstruction of the extrapolated variables used in the FOM to develop LADER scheme (see \cite{BFTVC17}).
From the computational point of view, after the calculus of  tensor $\boldsymbol{\varphi}_{j}\otimes \boldsymbol{\varphi}_{k}$ at the nodes, $N_{l}$, of the finite volume mesh, its divergence is approximated at each primal tetrahedra, $T_{m}$, by considering a Galerkin approach. Next, the value at each node is taken as the average of the values at the two tetrahedra from which it belongs. Finally, the $L^{2}$ product by the basis function, $\boldsymbol{\varphi}_{i}$, is performed. 

The storage of a third order tensor is only one of the possible choices in order to deal with the non-linear term of the Navier-Stokes equations. Such an approach is preferred when a relatively small number of basis functions is required, as in the cases examined in the present work. In cases with a large number of basis functions the storage of such a large (dense) tensor may lead to high storage costs and thus become unfeasible. In these cases it is possible to rely on alternative approaches for non-linear treatment such as the empirical interpolation method \cite{BARRAULT2004667}, the Gappy POD \cite{Everson1995} or the GNAT \cite{Carlberg2013623}.

% % % % % % % % % % % % % % % % % % % % % % % % % % % % % %
\subsubsection{Poisson equation} \hfill\\

In the framework of finite elements, it is well known that when using a mixed formulation, for the resolution of Navier-Stokes equations, to avoid pressure instabilities, the approximation spaces need to verify the inf-sup condition (see \cite{Boffi2013}). Also at the reduced order level, despite the snapshots are obtained by stable numerical methods, we can not guarantee that their properties will be preserved after the Galerkin projection and, with a saddle point formulation, one has to fulfill a reduced version of the inf-sup condition \cite{Rozza2007,BMQR15,SR17}. Instabilities at the reduced order level, and especially pressure instabilities, have been addressed by several authors, we recall here \cite{NPP05,CIJS14598,BBI09516,SR17}.

In this work, to be also consistent with the procedure used in the full order solver, instead of relying on a saddle point structure, we will substitute the divergence free condition by a Poisson equation for pressure following the works on \cite{OID86,JL04,ANR09,SHMLR17,SR17}. 

Taking the divergence of momentum equation, \eqref{eq:momentos}, and applying the divergence free condition yields the Poisson pressure equation
\begin{equation}
\Delta \pi = - \mathrm{div}\left[ \mathrm{div}\,  \left(\mathbf{w_{u}} \otimes \mathbf{w_{u}} \right)\right]  + \mathrm{div} \, \mathbf{f}_{\mathbf{u}}.\label{eq:poisson}
\end{equation}
The boundary conditions are determined by enforcing the divergence free constraint at the boundary
\begin{equation}
\frac{\partial \pi}{\partial \boldsymbol{\eta}} = -\mu \, \mathrm{rot}\left(\mathrm{rot} \mathbf{w}_{\mathbf{u}}\right) \cdot \boldsymbol{\eta} - \mathbf{g}_{t} \cdot \boldsymbol{\eta}
\label{eq:poisson_bc}
\end{equation}
with
\begin{equation}\mathbf{g}(\mathbf{x},t)= \mathbf{w_u}(\mathbf{x},t) \textrm{ in } \Gamma.\end{equation}
Alternative ways to enforce the boundary conditions can be found, for instance, in \cite{GS87}, \cite{Guer06} and \cite{LLP10}.

Substituting \eqref{eq:lc_velocity}, \eqref{eq:lc_source} and \eqref{eq:lc_pressure}  into equations \eqref{eq:poisson}-\eqref{eq:poisson_bc}, projecting the pressure equation onto the subspace spanned by the pressure modes, $X^{\mathrm{POD}_{\pi}}_{N_{\pi}}$, 
and taking into account the boundary conditions within the integration by parts,
we obtain the reduced system for the pressure
\begin{equation}
\mathbf{N}\mathbf{b}=\mathbf{a}^{T} \mathbf{D}\mathbf{a}  + \mathbf{H}\mathbf{a}+ \mathbf{P}\mathbf{a}+\mathbf{G}\mathbf{a},
\label{eq:dyn_sys_pi} 
\end{equation}
where
\begin{gather}
N_{ij}= \left< \mathrm{grad}\, \psi_{i}, \mathrm{grad}\, \psi_{j} \right>_{L^{2}},\quad
D_{ijk}= \langle\!\langle \psi_{i}, \frac{1}{\rho}\mathrm{div}\, \left[ \mathrm{div}\,\left(  \boldsymbol{\varphi}_{j}\otimes \boldsymbol{\varphi}_{k}\right) \right]  \rangle\!\rangle_{L^{2}}, \notag \\
H_{ij}=  \langle\!\langle \psi_{i}, \mathrm{div}\, \boldsymbol{\varsigma}_{j} \rangle\!\rangle_{L^{2}},\quad
P_{ij}=  \langle \boldsymbol{\eta}\times \psi_{i}, \frac{\mu}{\rho} \,\mathrm{rot}\,  \boldsymbol{\varphi}_{j}  \rangle_{\Gamma},\quad
G_{ij}=  \langle \psi_{i}, \mathbf{g}_{t} \cdot \boldsymbol{\eta} \rangle_{\Gamma}.
\end{gather}

Despite \eqref{eq:dyn_sys_w}  and \eqref{eq:dyn_sys_pi} are initially coupled, assuming 
that $\mathbf{N}$ is not singular, we get
\begin{equation}
\mathbf{b}=\mathbf{N}^{-1}\left( \mathbf{a}^{T} \mathbf{D}\mathbf{a}  + \mathbf{H}\mathbf{a}+ \mathbf{P}\mathbf{a}+\mathbf{G}\mathbf{a} \right) . \label{eq:computation_pressurecoef}
\end{equation}
Therefore, substituting \eqref{eq:computation_pressurecoef} on \eqref{eq:dyn_sys_w}, it results
\begin{equation}
\dot{\mathbf{a}} = -\mathbf{a}^{T} \mathbf{C}\mathbf{a}+ \mathbf{B} \mathbf{a} - \mathbf{K}\, \mathbf{N}^{-1}\left( \mathbf{a}^{T} \mathbf{D}\mathbf{a}  + \mathbf{H}\mathbf{a}+ \mathbf{P}\mathbf{a}+\mathbf{G}\mathbf{a} \right)   + \mathbf{F} \mathbf{a}.\label{eq:dyn_sys_wpi}
\end{equation}
So, the vectors of coefficients can be computed in two stages.
Firstly, we obtain the coefficients of the velocity, $\mathbf{a}$, by solving the algebraic system \eqref{eq:dyn_sys_wpi}. Next, the coefficients related to the pressure, $\mathbf{b}$, are computed from  \eqref{eq:computation_pressurecoef}.

\begin{remark}
	It is important to notice that the solution obtained for the pressure, as in the FOM solution, since we are working with the gradient of the pressure instead with the pressure itself, is defined up to an arbitrary constant. To correct this issue it is necessary to impose an initial condition for the pressure. As initial condition for pressure we use the projection of the FOM initial pressure solution onto the POD spaces. 
	The enforcement of an initial condition for pressure permits to obtain the constant and therefore to ensure the consistency of the FOM and the ROM.
\end{remark}

% % % % % % % % % % % % % % % % % % % % % % % % % % % % % %
\subsubsection{Transport equation}\hfill\\

Performing Galerkin projection onto the transport equation gives
\begin{equation}
\dot{\mathbf{c}} + \mathbf{a}^{T} \mathbf{E} \mathbf{c} -\mathbf{Q} \mathbf{c} =0 \label{eq:dyn_sys_y} 
\end{equation}
where
\begin{gather}
\mathbf{E}_{ijk}= \langle \chi_{i}, \frac{1}{\rho}\mathrm{div}\left( \boldsymbol{\varphi}_{j}\, \chi_{k} \right)  \rangle_{L^{2}},\quad
\mathbf{Q}_{ij}=  \langle \chi_{i}, \mathcal{D}\Delta \chi_{j} \rangle_{L^{2}}.
\end{gather}
The former algebraic system can be solved in a coupled way with \eqref{eq:dyn_sys_wpi} so that the vectors of coefficients $\mathbf{a}$ and $\mathbf{c}$ are obtained within the same stage.

% % % % % % % % % % % % % % % % % % % % % % % % % % % % % %
% % % % % % % % % % % % % % % % % % % % % % % % % % % % % %
\subsection{Lifting function}\hfill\\

So far, we have assumed homogeneous boundary conditions. In case non-homogeneous Dirichlet boundary conditions are imposed, we can homogenize the original snapshots by defining a lifting function. 

The lifting function method firstly proposed in \cite{Graham1999} for boundary condition that can be parametrized by a single multiplicative coefficient, and generalized in \cite{Gunzburger2007} for generic functions, is only one of the possible methods to consider boundary conditions also at the reduced order level. Some of other possible approaches consists in the modified basis function method \cite{Gunzburger2007} or in the penalty method approach \cite{Sirisup2005}. In the follows the lifting function method is briefly recalled. 

Let us assume that we have a constant on time non zero Dirichlet boundary condition for the velocity on the boundary. Then, a possible election would be to define the lifting function as the mean function of the velocity along the snapshots,
\begin{equation}
\boldsymbol{\varphi}_{\mathrm{lift}}\left(\mathbf{x}\right)=\frac{1}{N_{s}} \sum_{n=1}^{N_{s}} \mathbf{w}_{\mathbf{u}}^{n}\left(\mathbf{x}\right).
\end{equation}
Next, the homogeneous snapshots are computed as
\begin{equation}
\hat{\mathbf{w}}_{\mathbf{u}}^{n}\left(\mathbf{x}\right)= \mathbf{w}_{\mathbf{u}}^{n}\left(\mathbf{x}\right)-
\boldsymbol{\varphi}_{\mathrm{lift}}\left(\mathbf{x}\right), \quad n=1,\dots,N_{s}.
\end{equation}
Following the methodology already presented in Section \ref{sec:pod_spaces}, we obtain the POD space related to the homogenized snapshots,
\begin{equation}
\hat{X}^{\mathrm{POD}}_{N-1} := \displaystyle{\textrm{span}_{i=1,\dots,N-1} } \left\lbrace\hat{\boldsymbol{\varphi}}_i\right\rbrace.
\end{equation}
Therefore, the final POD space for the original set of snapshots is given by
\begin{equation}
X^{\mathrm{POD}}_{N} = 
\displaystyle{\textrm{span} }\left\lbrace \boldsymbol{\varphi}_{\mathrm{lift}}, \left\lbrace  \hat{\boldsymbol{\varphi}}_i, \, i=1,\dots,N-1\right\rbrace\right\rbrace.
\end{equation}

Since the lifting function may not be orthogonal to the remaining elements of the basis, the mass matrix of the system related to the momentum equation must be computed. Accordingly, the modified algebraic system to be solved reads
\begin{equation}
\mathbf{M}\dot{\mathbf{a}} = -\mathbf{a}^{T} \mathbf{C}\mathbf{a}+ \mathbf{B} \mathbf{a} - \mathbf{K} \mathbf{b}  + \mathbf{F} \mathbf{a},\label{eq:dyn_sys_w_lifting}
\end{equation}
with
\begin{equation}
\mathbf{M}_{ij}=\langle \boldsymbol{\varphi}_{i},\boldsymbol{\varphi}_{j}\rangle_{L^{2}} .
\end{equation}

% % % % % % % % % % % % % % % % % % % % % % % % % % % % % %
% % % % % % % % % % % % % % % % % % % % % % % % % % % % % %
\subsection{Initial conditions}\hfill\\

The initial conditions for the ROM systems of ODEs \eqref{eq:dyn_sys_w_lifting} and  \eqref{eq:dyn_sys_y} are computed by performing a Galerkin projection of the initial solution, $\mathbf{w}\left(\mathbf{x},t^{1}\right)$, onto the POD basis:
\begin{equation}
a_{0\, i}=\langle \boldsymbol{\varphi}_{i}\left(\mathbf{x}\right) ,\mathbf{w}_{\mathbf{u}}\left(\mathbf{x},t^{1}\right)
\rangle_{L^{2}},\quad
c_{0\, i}=\langle \chi_{i}\left(\mathbf{x}\right) , w_{\mathrm{y}}\left(\mathbf{x},t^{1}\right)\rangle_{L^{2}}.
\end{equation}
If non-orthogonal basis functions are considered for the linear momentum, the initial coefficients, $\mathbf{a}_{0}$, must be obtained solving the linear system of equations
\begin{equation}
\mathbf{M} \mathbf{a}_{0} =\mathbf{e}
\end{equation}
with $e_{i}=\langle \boldsymbol{\varphi}_{i}\left(\mathbf{x}\right) ,\mathbf{w}\left(\mathbf{x},t^{1}\right) \rangle_{L^{2}}$.

Similarly, the components of the initial pressure coefficient, $\mathbf{b}_{0}$, result
\begin{equation}b_{0\,i}=\langle \psi_{i},\pi\left(\mathbf{x},t^{1}\right) \rangle_{L^{2}}.
\end{equation}

\vspace{0.2cm}
% % % % % % % % % % % % % % % % % % % % % % % % % % % % % %
%                  Numerical results                      %
% % % % % % % % % % % % % % % % % % % % % % % % % % % % % %
\section{Numerical results}\label{sec:numerical_results}

In this section we present the results obtained for two different test problems. In both numerical tests we simulate the same conditions for the FOM and the ROM. Therefore, the ROM solutions are compared against the high fidelity ones.

\subsection{Manufactured test}\hfill\\

The first test to be considered has been obtained by using the method of the manufactured solutions. We consider the computational domain $\Omega=\left[0,1\right]^{3}$ and the flow being defined by
\begin{gather}
\rho=1, \quad \pi= t \cos \left( \pi \left(x+y+z\right) \right),\quad
\mathbf{w}_{\mathbf{u}}= \left(\cos^{2}\left(\pi x t\right),
e^{-2\pi y t},-\cos\left(\pi x t\right)\right)^{T},
\end{gather} 
with $\mu=10^{-2}$ and the source terms given by
\begin{gather}
f_{u_1}= \pi y  \cos(\pi t y)  \cos(\pi t z) 
- \pi z  \sin(\pi t y)  \sin(\pi t z) 
+ 2  \pi^2 t^2 \mu  \sin(\pi t y)  \cos(\pi t z)\\ 
- \pi t  \cos(\pi t y)  \cos^2(\pi t z) 
- \pi t  \sin(\pi t y)  \sin(\pi t z) e^{-2  \pi t^2 x} 
- t \pi  \sin(\pi (x + y + z)) ,\notag\\
f_{u_2}= \pi z  \sin(\pi t z) 
- \pi^2 t^2 \mu  \cos(\pi t z) 
+ \pi t  \sin(\pi t z) e^{-2  \pi t^2 x}
- t \pi  \sin(\pi (x + y + z)),\\
f_{u_3}= - 4  \pi^2 t^4 \mu e^{-2  \pi t^2 x}
- 4  \pi t x e^{-2  \pi t^2 x}- t \pi  \sin(\pi (x + y + z))\\
- 2  \pi t^2  \sin(\pi t y) e^{-2  \pi t^2 x } \cos(\pi t z).\notag
\end{gather}

\begin{figure}[h]
	\centering
	\includegraphics[width=0.4\linewidth]{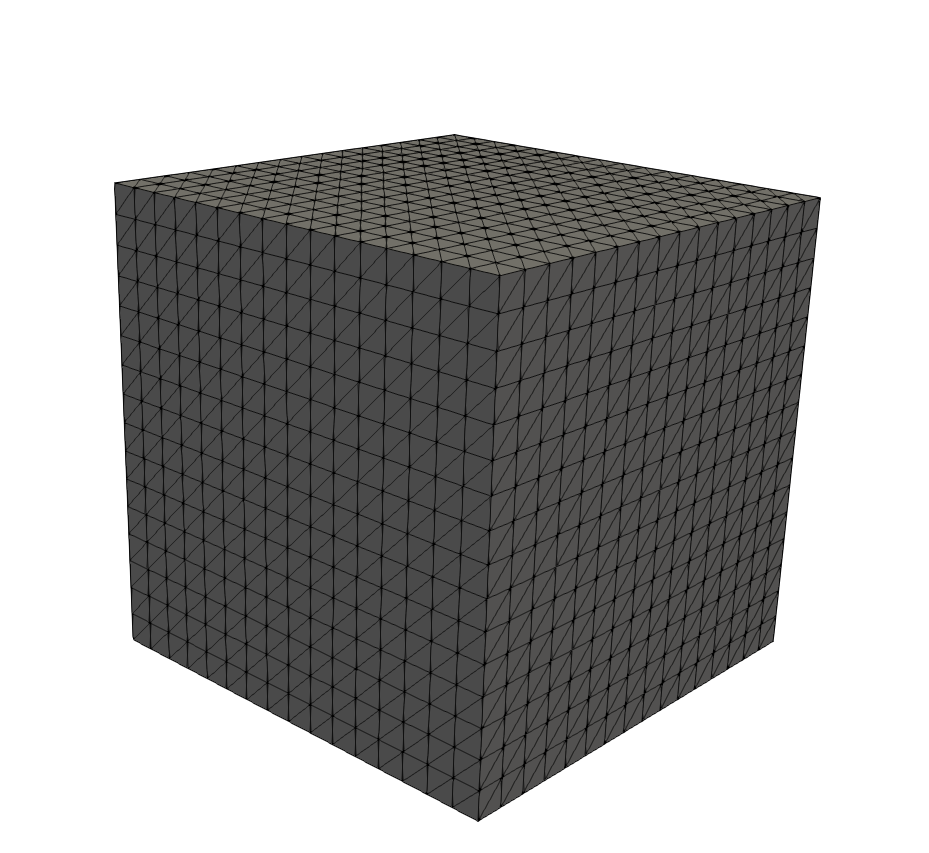}
	\caption{Manufactured test. Primal tetrahedral mesh.}
	\label{fig:mallacubo16ef}
\end{figure}

We generate an initial tetrahedral mesh of $24576$ elements
whose dual mesh accounts for $50688$ nodes (see Figure \ref{fig:mallacubo16ef}).
The minimum and maximum volumes of the dual cells are $1.02E-05$ and $2.03E-05$, respectively. The time interval of the simulation is taken to be $T=\left[0,2.5\right]$.

The full order model simulation was run considering the second order in space and time LADER methodology (\cite{TMN01}, \cite{BTVC16}, \cite{BFTVC17}). 
We consider a $CFL=5$ to determine the time step at each iteration (see \cite{BFTVC17} for further details on its computation).

\begin{table}[h]
	\centering
	\renewcommand{\arraystretch}{1.3}
	\begin{tabular}{|c||c|c|}
		\hline 
		Number of modes & $\mathbf{w}_{\mathbf{u}}$& $\pi$ \\ \hline
		\hline 
		$1$& $0.4232018 $& $\mathbf{0.9999753}$ \\  
		$2$& $0.7194136 $& $0.9999924$ \\ 
		$3$& $0.9276722 $& $0.9999959$ \\ 
		$4$& $0.9724632 $& $0.9999975$ \\ 
		$5$& $0.9910768 $& $0.9999984$ \\ 
		$6$& $0.9989818 $& $0.9999988$ \\ 
		$7$& $0.9999167 $& $0.9999991$ \\ 
		$8$& $0.9999800 $& $0.9999994$ \\ 
		$9$& $\mathbf{0.9999955}$& $0.9999995$ \\ \hline 
	\end{tabular} 
	
	\vspace{0.2cm}
	\caption{Manufactured test. Cumulative eigenvalues for the linear momentum and the pressure. The values in which  the fixed bounds $\kappa_{\mathbf{w}_{\mathbf{u}}}\!=\!99.999\%$ and
		$\kappa_{\pi}\!=\!99.99\%$ are reached are written in bold font.}
	\label{tab:mms_eigen}
\end{table}

The snapshots are taken every $0.01$s given a total number of snapshots equal to $250$.
The dimension of the reduced spaces are $N=9$ and $N_{\pi}=1$. This selection is done to retain more than $\kappa_{\mathbf{w}_{\mathbf{u}}}\!=\!99.999\%$  of the energy and $\kappa_{\pi}\!=\!99.99\%$. Table \ref{tab:mms_eigen} contains the cumulative eigenvalues obtained.   

\begin{figure}[H]
	\centering
	\includegraphics[width=0.8\linewidth]{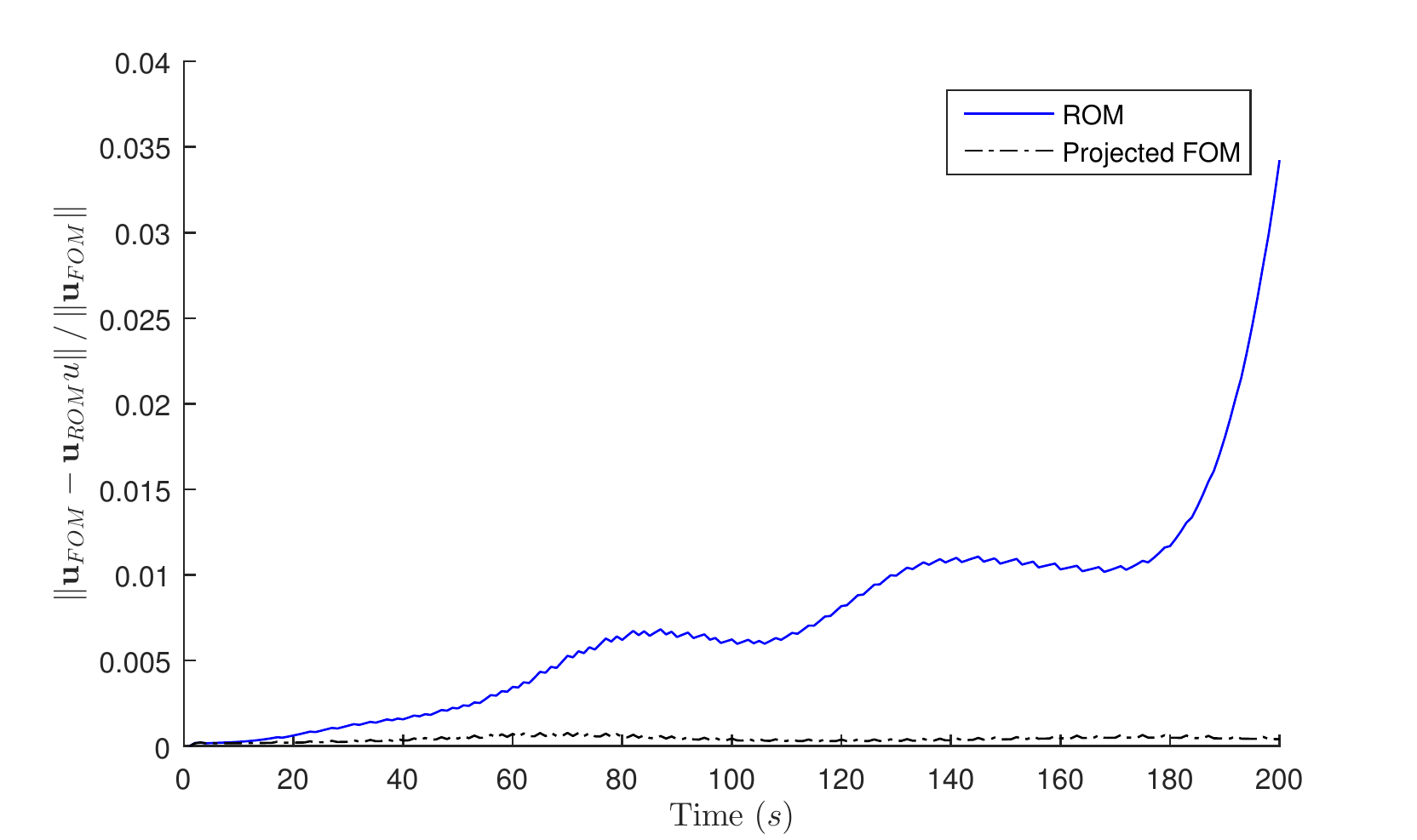}
	\caption{Manufactured test. Relative error of the linear momentum field for the projected FOM solution and the ROM solution. }
	\label{fig:wu_t4n_250}
\end{figure}

\begin{figure}[H]
	\centering
	\includegraphics[width=0.8\linewidth]{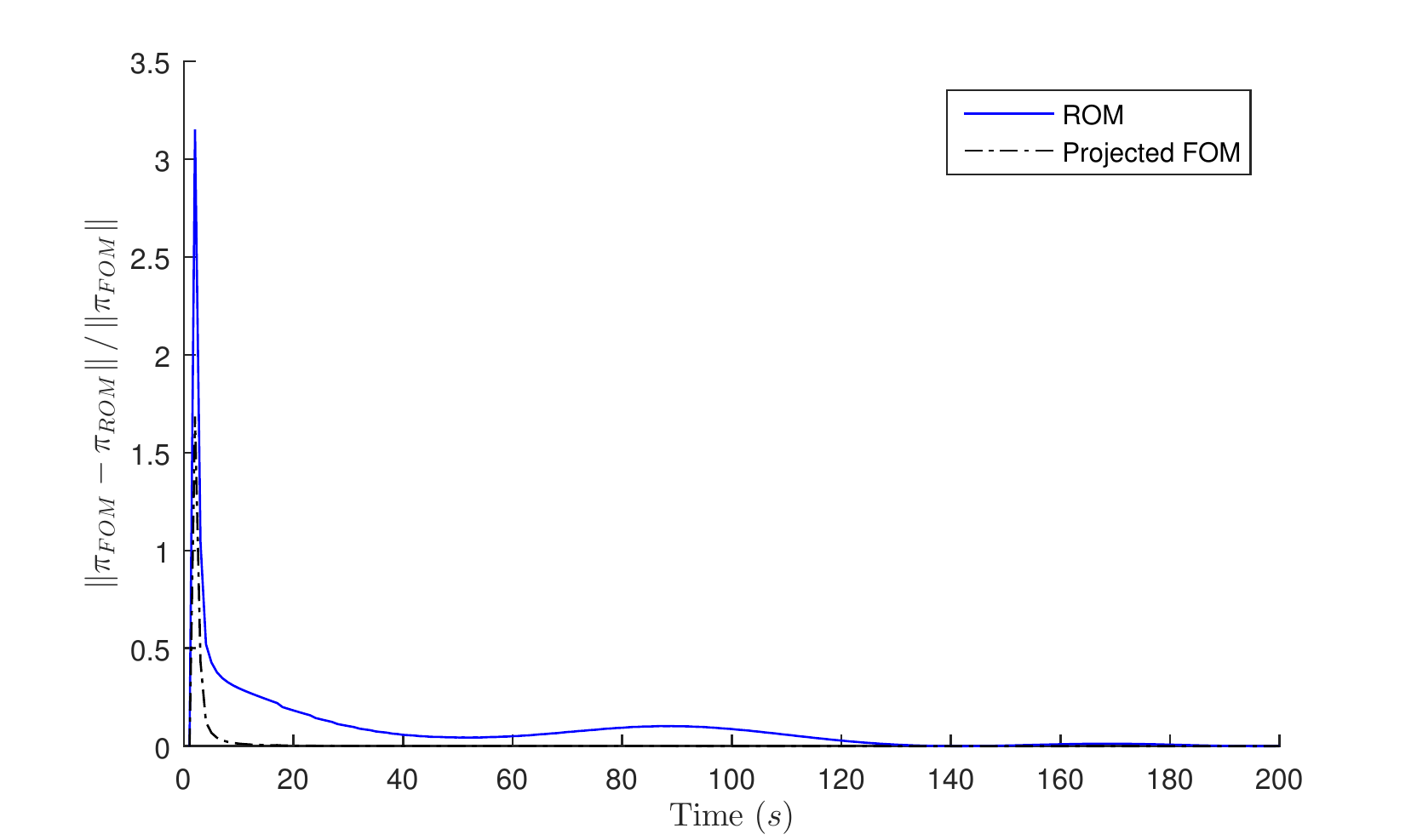}
	\caption{Manufactured test. Relative error of the pressure field for the projected FOM solution and the ROM solution. }
	\label{fig:pi_t4n_250}
\end{figure}

To assess the methodology the ROM solution is compared against the FOM solution. The $L^{2}$ error for the diverse snapshots is depicted in Figures \ref{fig:wu_t4n_250} and \ref{fig:pi_t4n_250}. Furthermore, the errors obtained by projecting the FOM solution onto the considered POD basis are also portrayed.
Figures \ref{fig:wumagnitudexplussl}, \ref{fig:wumagnitudeyplussl} and \ref{fig:wumagnitudezplussl}  show the solutions obtained using FOM and ROM for different time instants at the mid planes $x=0.5$, $y=0.5$ and $z=0.5$ respectively. A good agreement between the two solutions is observed.


\begin{center}
	\begin{minipage}{0.9\linewidth}
		\begin{figure}[H]
			\begin{center}
				{\small \hspace*{1cm}	$t=0.1s$\hfill $t=1s$ \hfill $t=1.5s$ \hfill $t=2s$ \hspace*{1.2cm}}
			\end{center}	
			\includegraphics[width=0.24\linewidth]{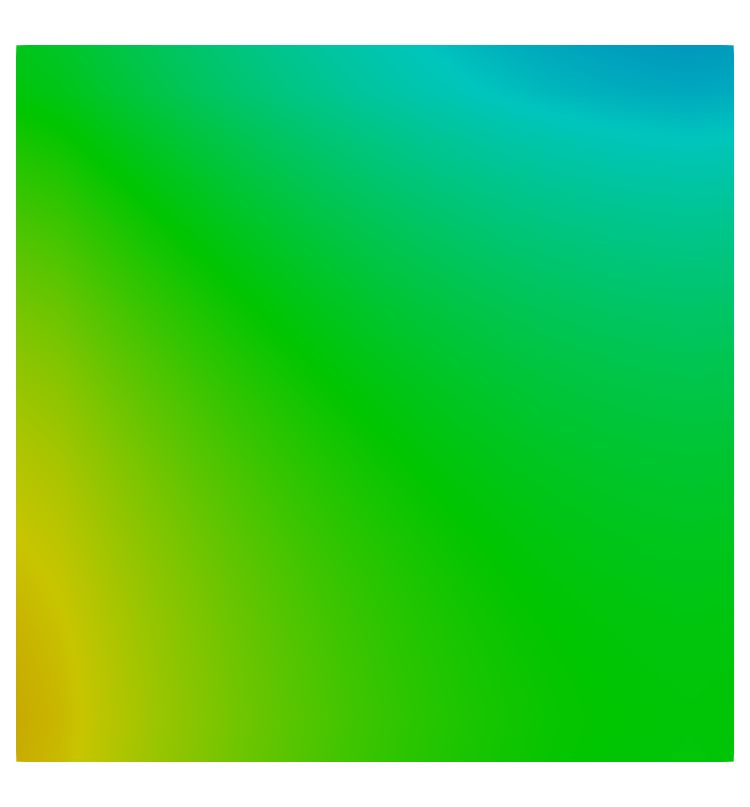}
			\hfill
			\includegraphics[width=0.24\linewidth]{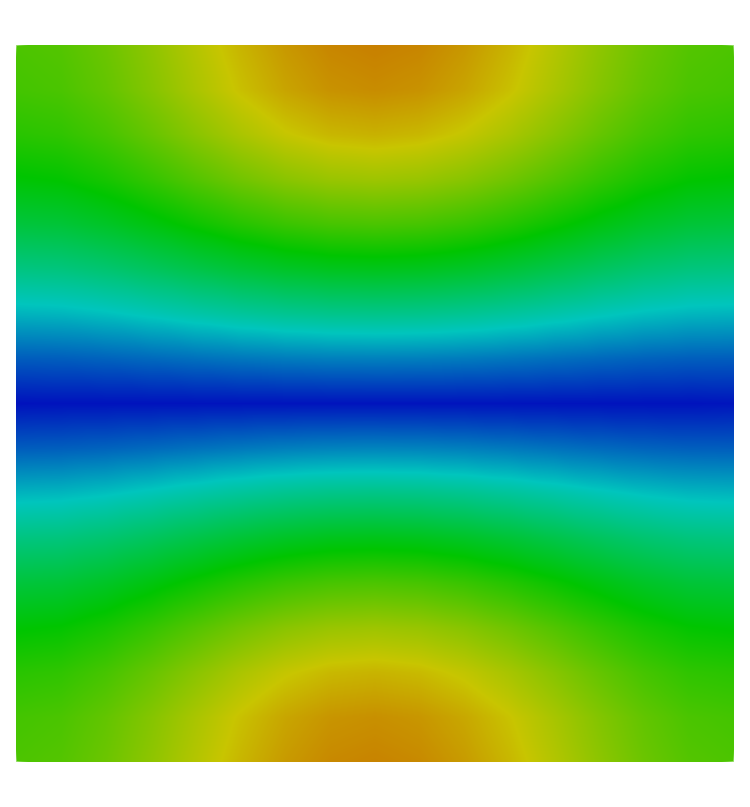}
			\hfill
			\includegraphics[width=0.24\linewidth]{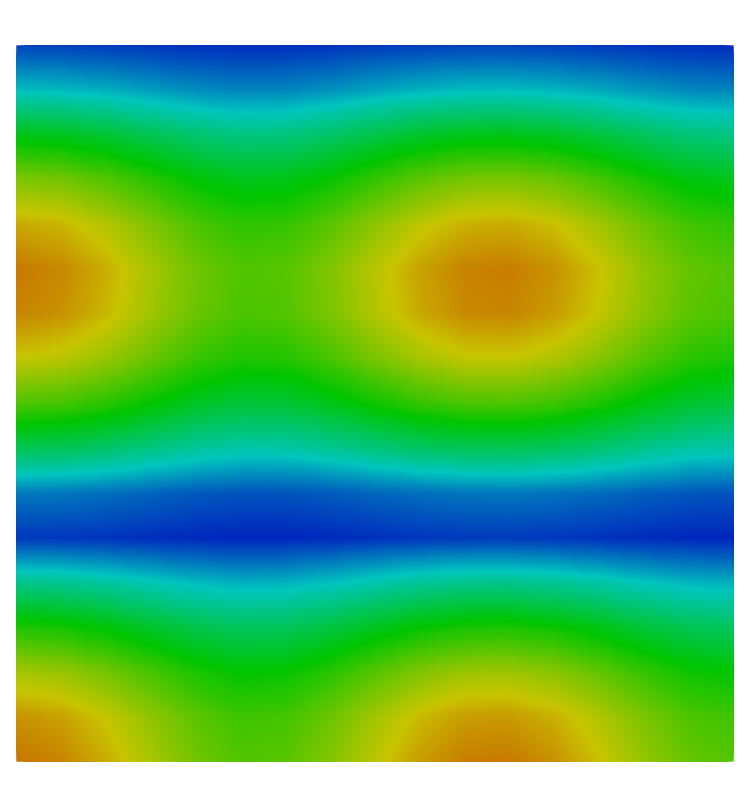}
			\hfill
			\includegraphics[width=0.24\linewidth]{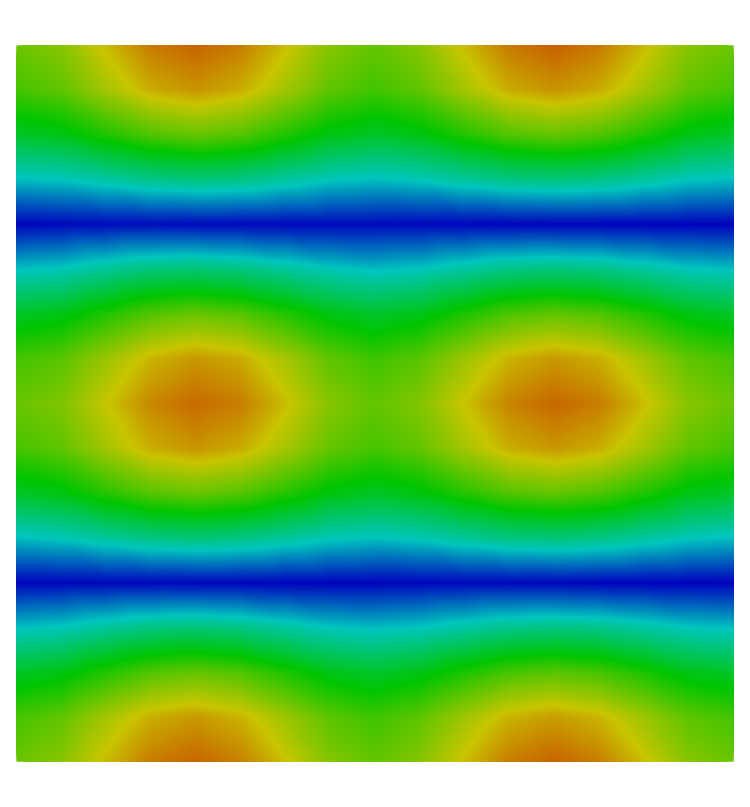}
			
			\includegraphics[width=0.24\linewidth]{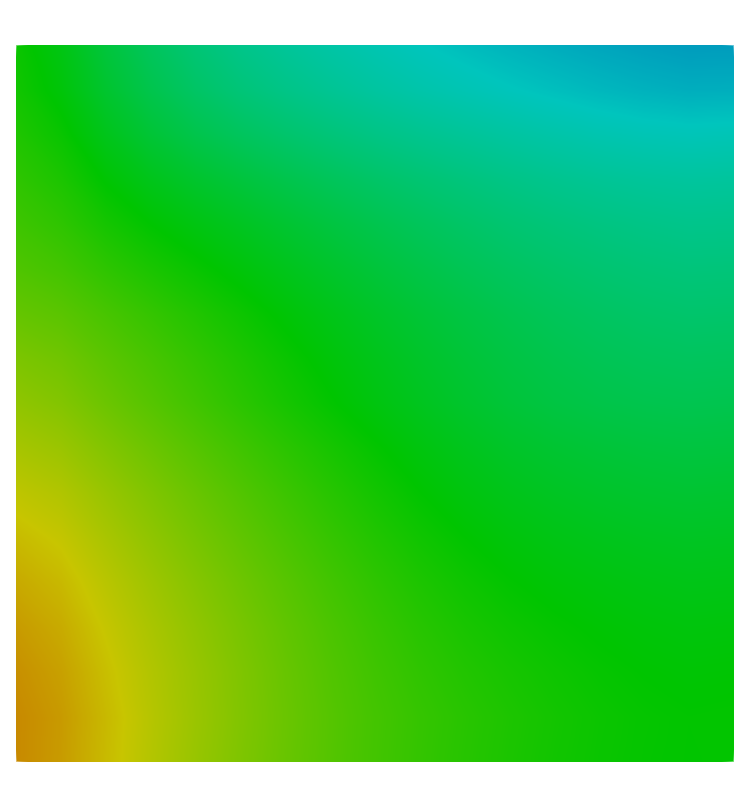}
			\hfill
			\includegraphics[width=0.24\linewidth]{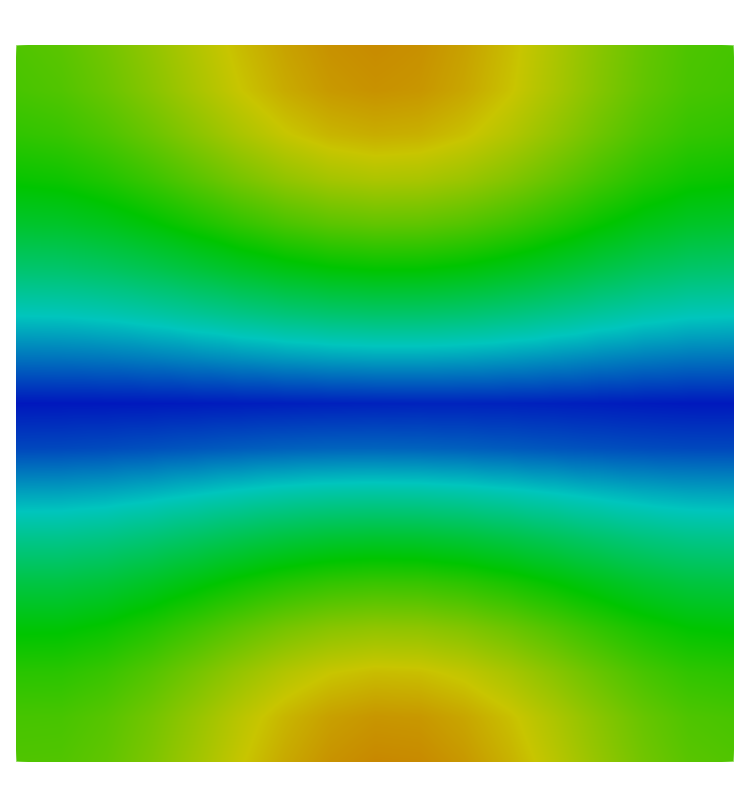}
			\hfill
			\includegraphics[width=0.24\linewidth]{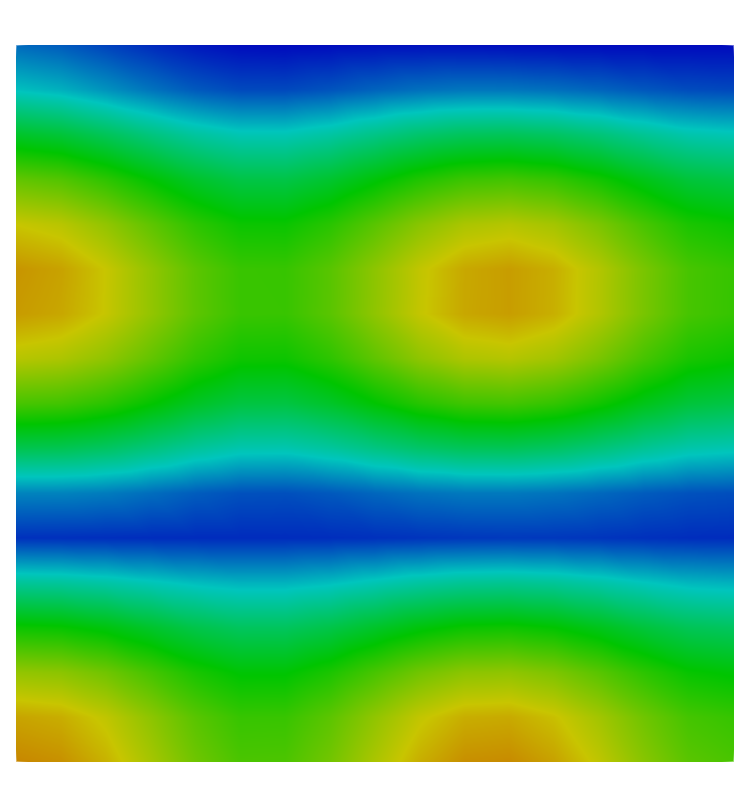}
			\hfill
			\includegraphics[width=0.24\linewidth]{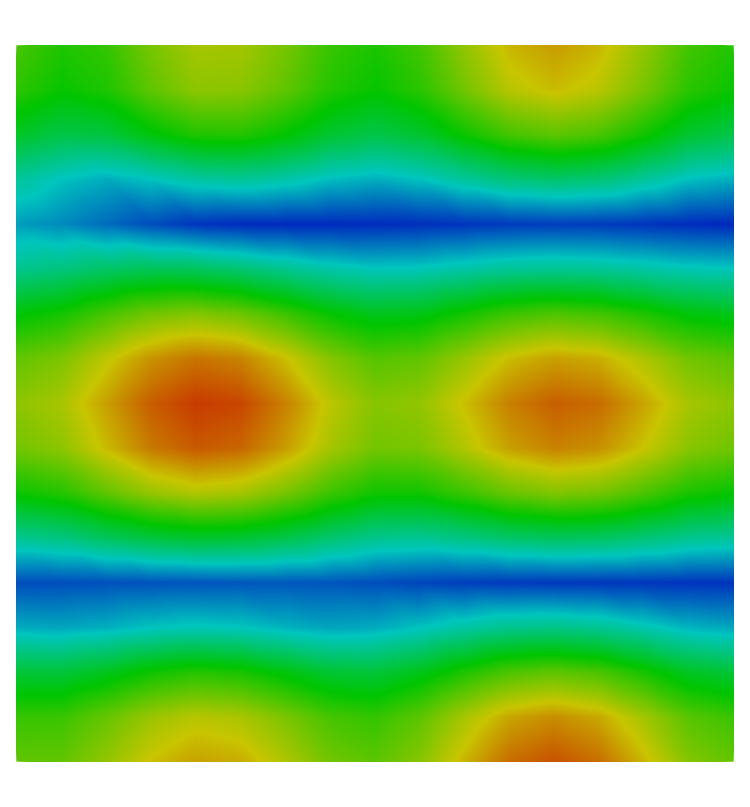}
			
			\vspace{0.2cm}
			\includegraphics[width=0.24\linewidth]{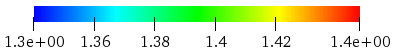}\hfill
			\includegraphics[width=0.24\linewidth]{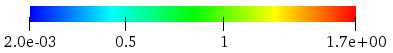}\hfill
			\includegraphics[width=0.24\linewidth]{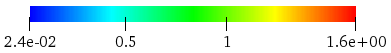}\hfill
			\includegraphics[width=0.24\linewidth]{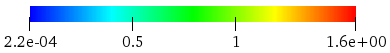}
			
			\caption{Manufactured test. Comparison of the velocity Magnitude at plane $x=0.5$.
				Top: FOM solution. Bottom: ROM approximation.}
			\label{fig:wumagnitudexplussl}
		\end{figure}
	\end{minipage}
\end{center}

\begin{center}
	\begin{minipage}{0.9\linewidth}
		\begin{figure}[H]
			\begin{center}
				{\small \hspace*{1cm}	$t=0.1s$\hfill $t=1s$ \hfill $t=1.5s$ \hfill $t=2s$ \hspace*{1.2cm}}
			\end{center}	
			\includegraphics[width=0.24\linewidth]{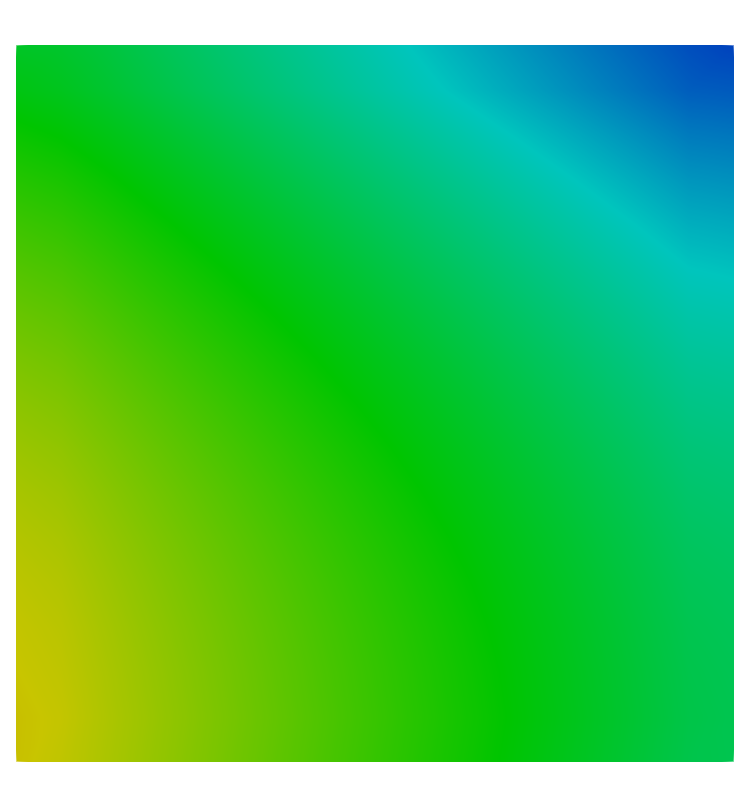}
			\hfill
			\includegraphics[width=0.24\linewidth]{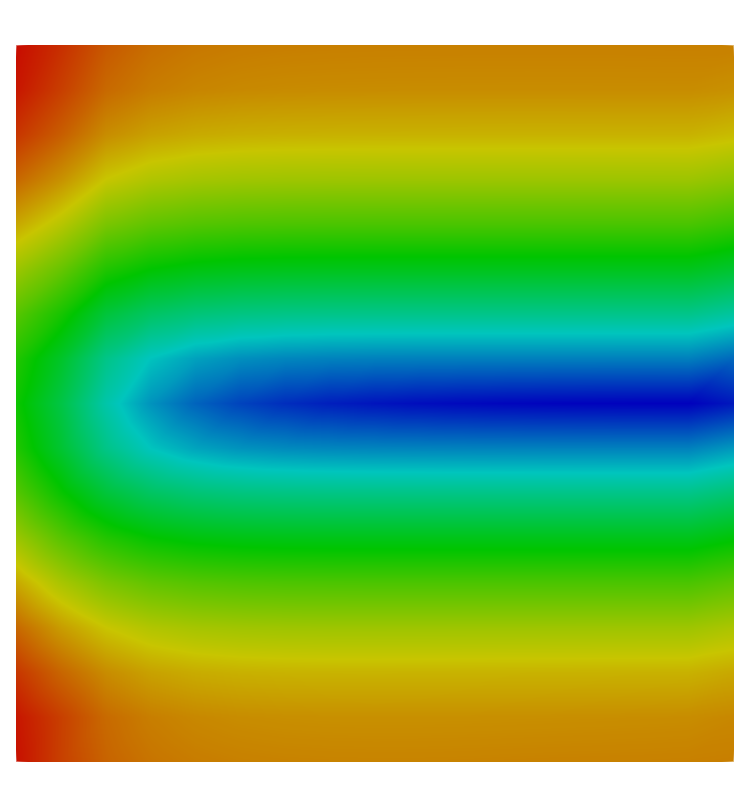}
			\hfill
			\includegraphics[width=0.24\linewidth]{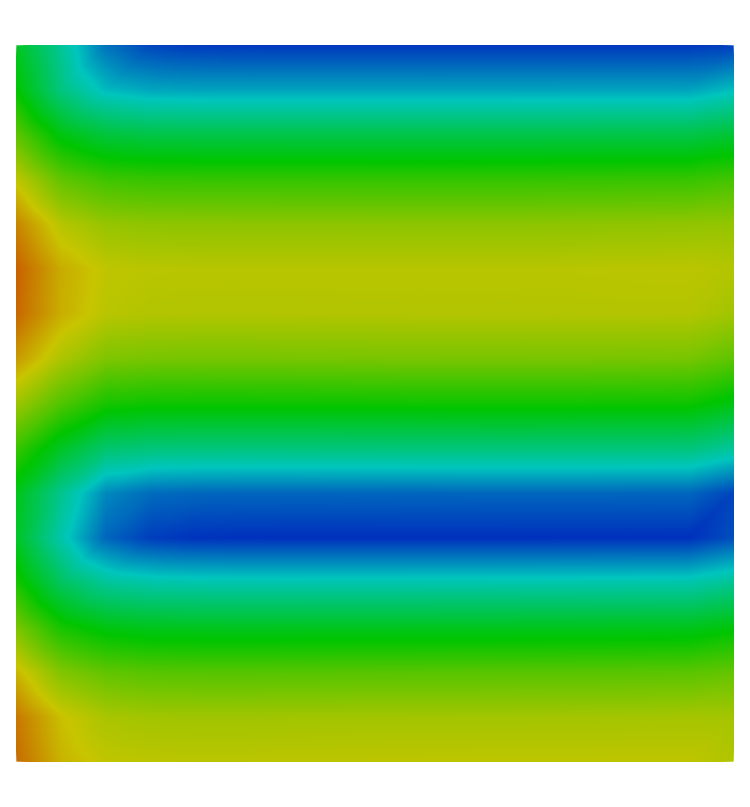}
			\hfill
			\includegraphics[width=0.24\linewidth]{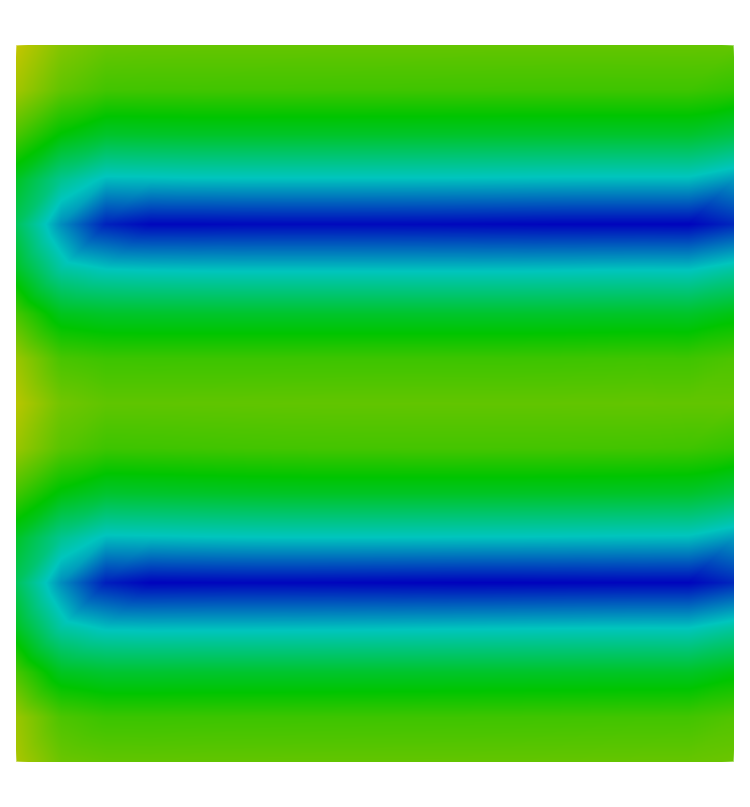}
			
			\includegraphics[width=0.24\linewidth]{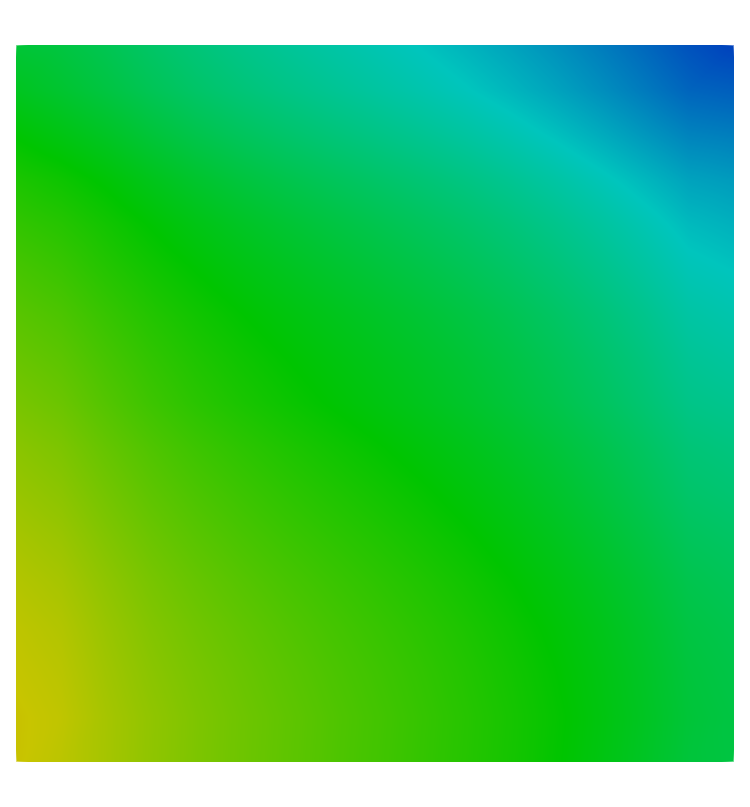}
			\hfill
			\includegraphics[width=0.24\linewidth]{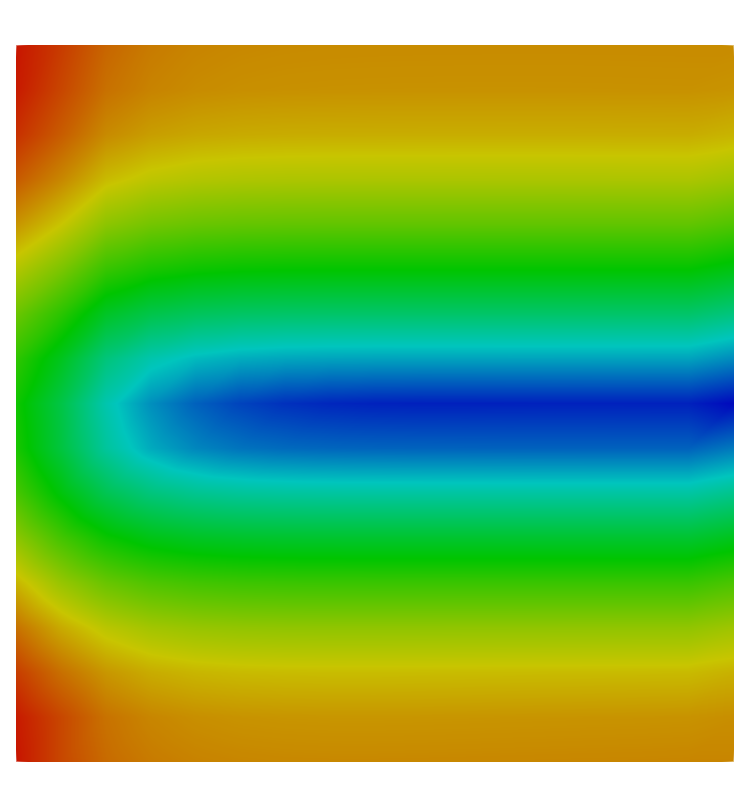}
			\hfill
			\includegraphics[width=0.24\linewidth]{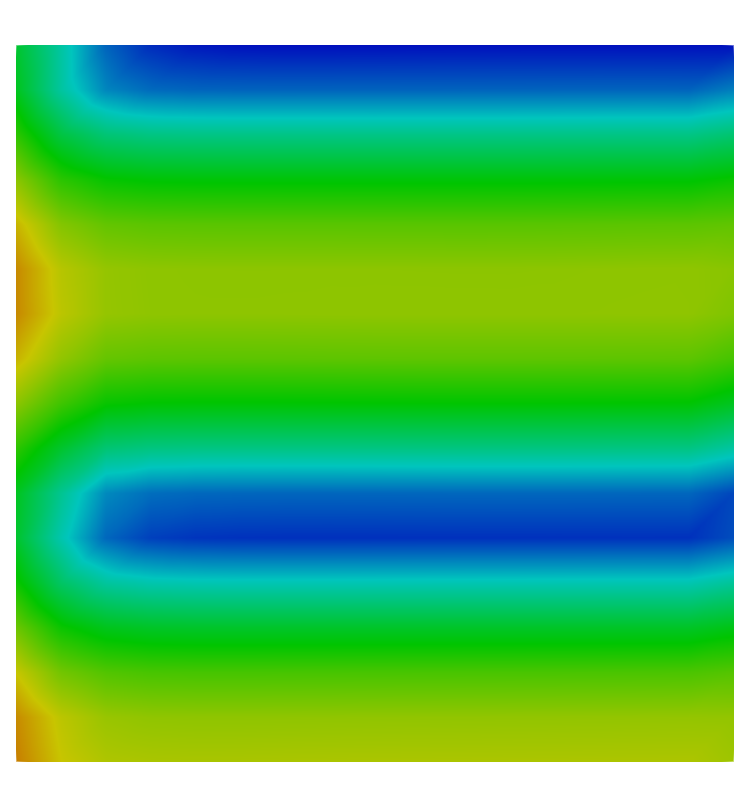}
			\hfill
			\includegraphics[width=0.24\linewidth]{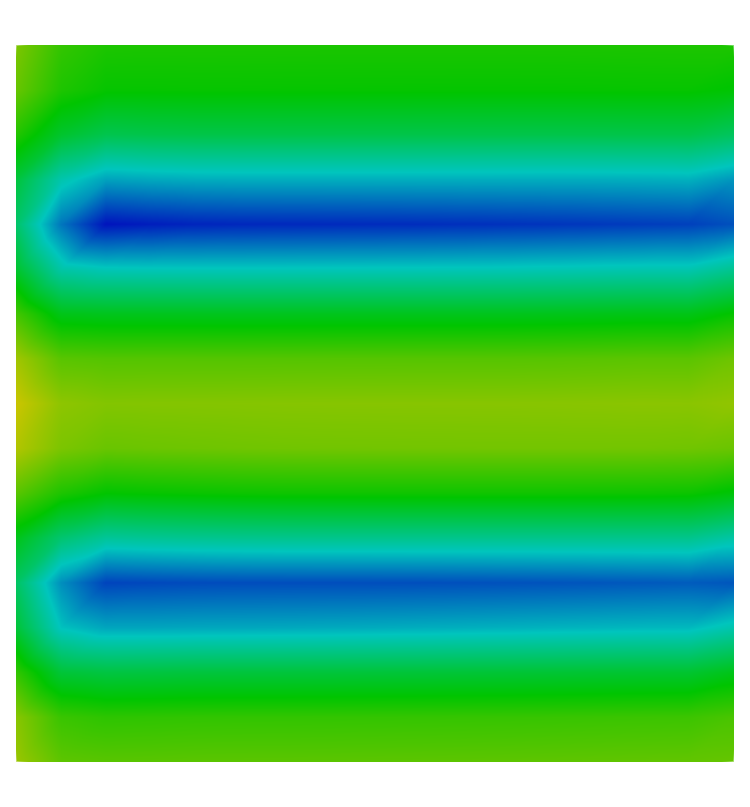}
			
			\vspace{0.2cm}
			\includegraphics[width=0.24\linewidth]{./legend001_sn}\hfill
			\includegraphics[width=0.24\linewidth]{./legend010_sn}\hfill
			\includegraphics[width=0.24\linewidth]{./legend015_sn}\hfill
			\includegraphics[width=0.24\linewidth]{./legend020_sn}

			\caption{Manufactured test. Comparison of the velocity Magnitude at plane $y=0.5$.
				Top: FOM solution. Bottom: ROM approximation.}
			\label{fig:wumagnitudeyplussl}
		\end{figure}
	\end{minipage}
\end{center}

\begin{center}
	\begin{minipage}{0.9\linewidth}
		\begin{figure}[H]
			\begin{center}
				{\small \hspace*{1cm}	$t=0.1s$\hfill $t=1s$ \hfill $t=1.5s$ \hfill $t=2s$ \hspace*{1.2cm}}
			\end{center}	
			\includegraphics[width=0.24\linewidth]{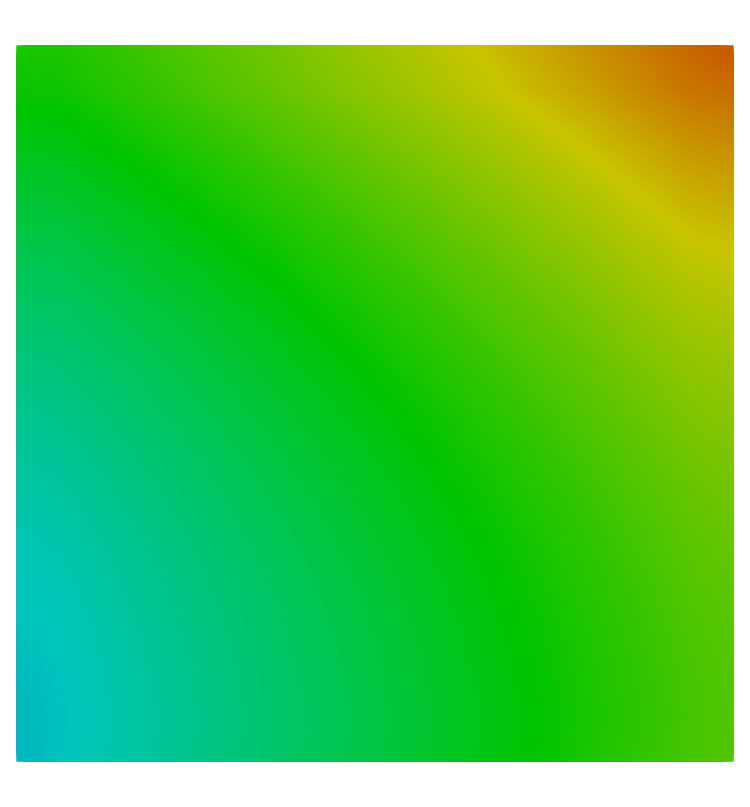}
			\hfill
			\includegraphics[width=0.24\linewidth]{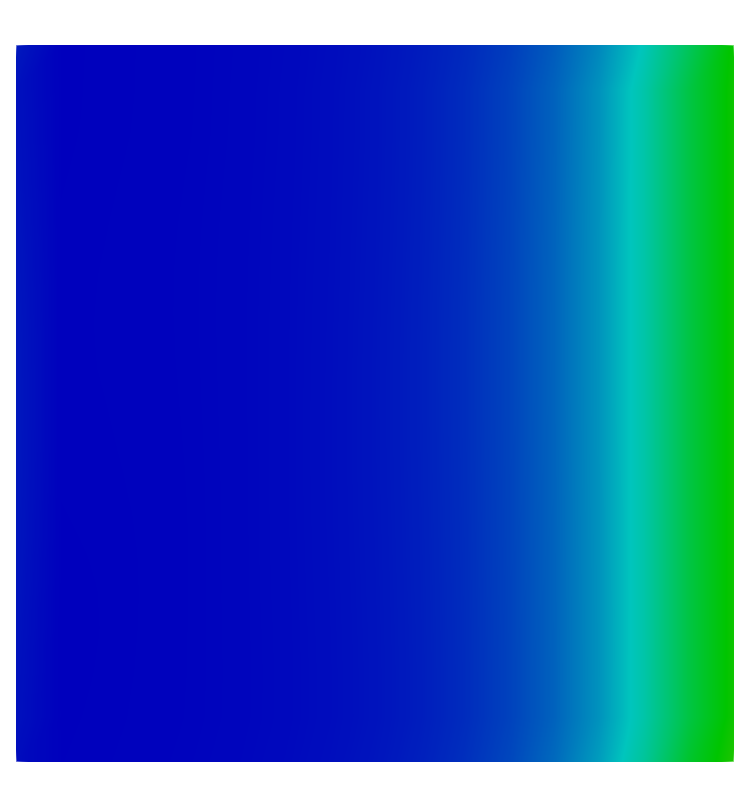}
			\hfill
			\includegraphics[width=0.24\linewidth]{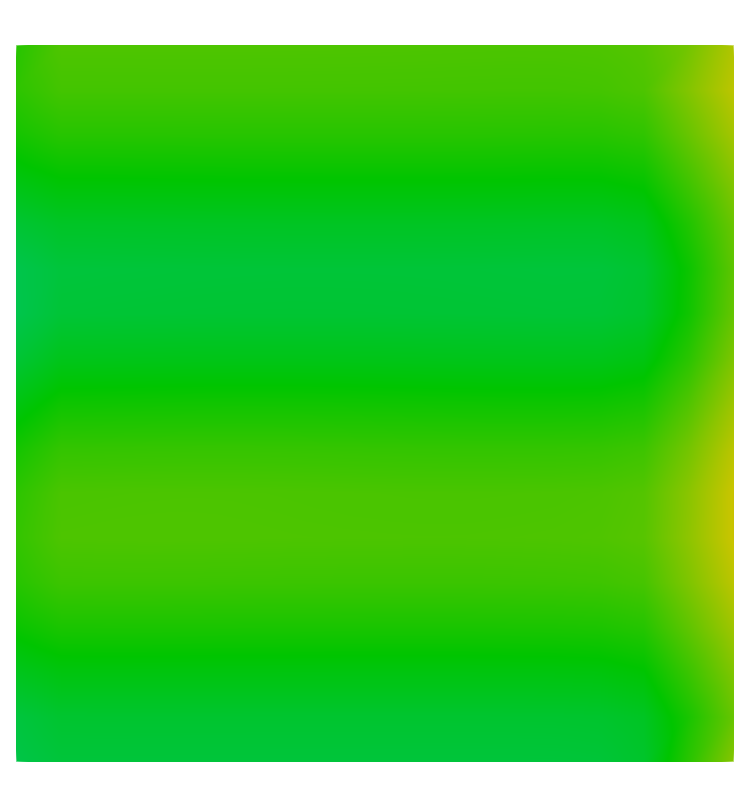}
			\hfill
			\includegraphics[width=0.24\linewidth]{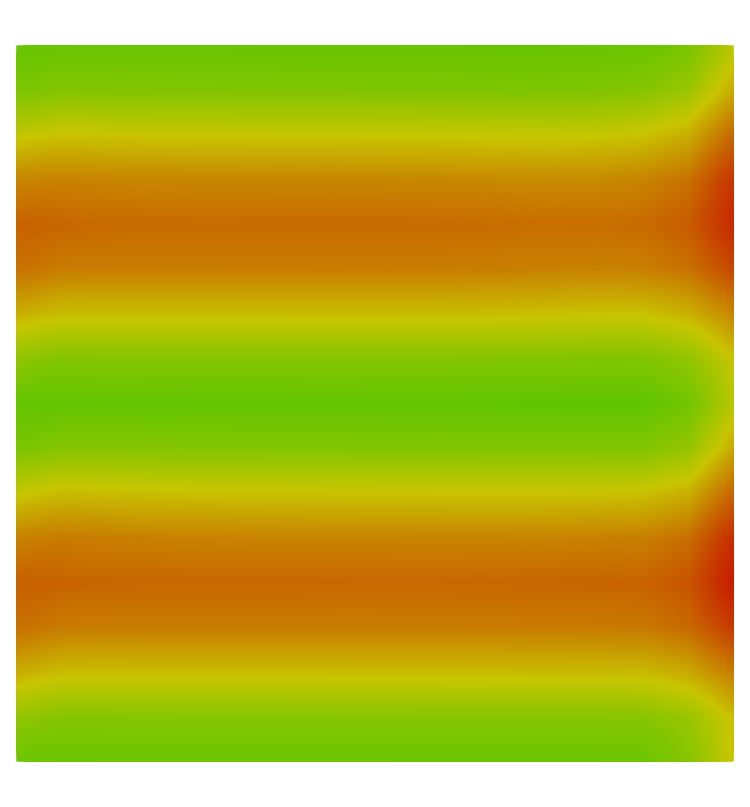}
			
			\includegraphics[width=0.24\linewidth]{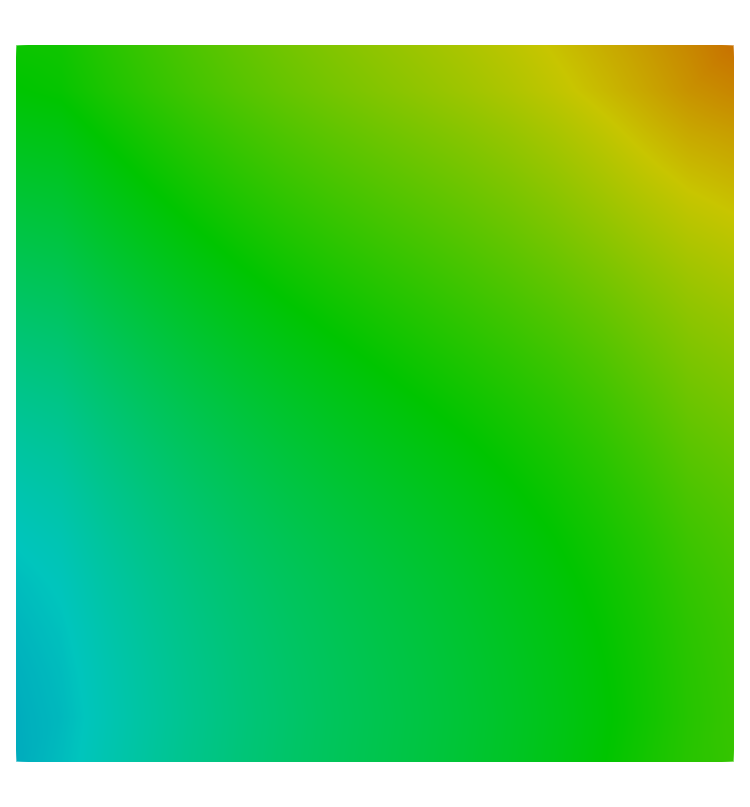}
			\hfill
			\includegraphics[width=0.24\linewidth]{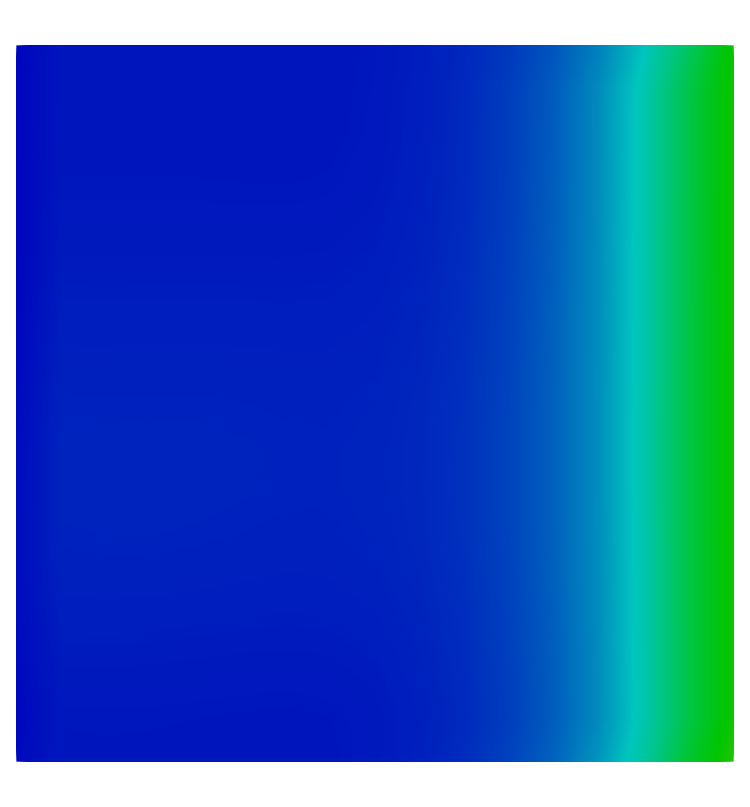}
			\hfill
			\includegraphics[width=0.24\linewidth]{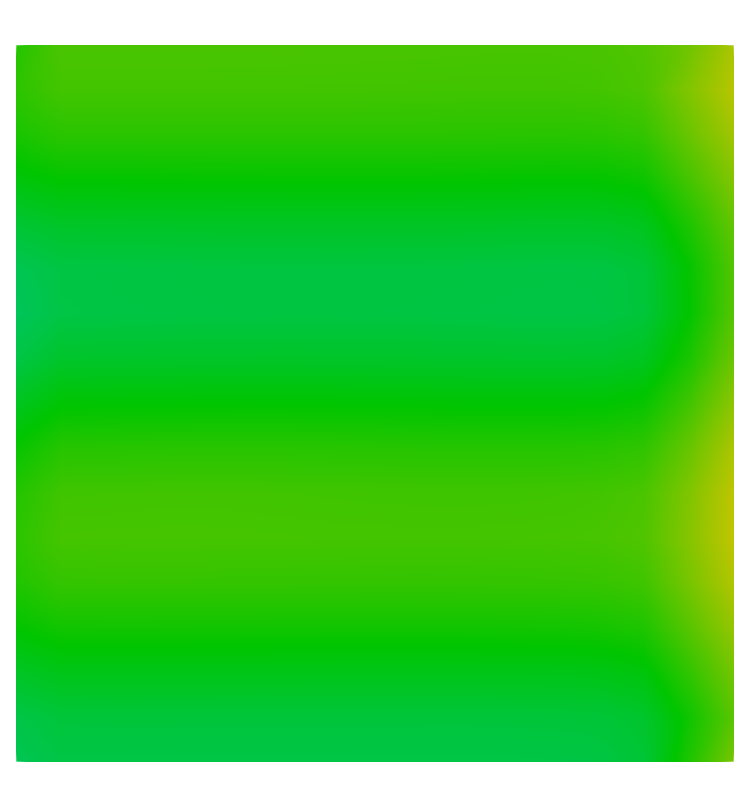}
			\hfill
			\includegraphics[width=0.24\linewidth]{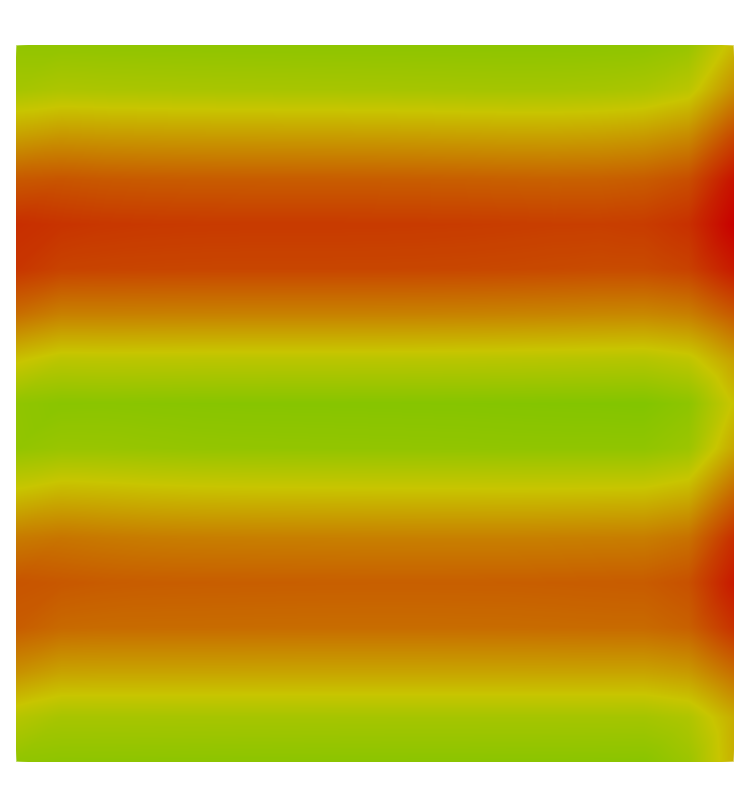}
			
			\vspace{0.2cm}
			\includegraphics[width=0.24\linewidth]{./legend001_sn}\hfill
			\includegraphics[width=0.24\linewidth]{./legend010_sn}\hfill
			\includegraphics[width=0.24\linewidth]{./legend015_sn}\hfill
			\includegraphics[width=0.24\linewidth]{./legend020_sn}
			
			\caption{Manufactured test. Comparison of the velocity Magnitude at plane $z=0.5$.
				Top: FOM solution. Bottom: ROM approximation.}
			\label{fig:wumagnitudezplussl}
		\end{figure}
	\end{minipage}
\end{center}

\begin{figure}[H]
	\centering
	\includegraphics[width=0.8\linewidth]{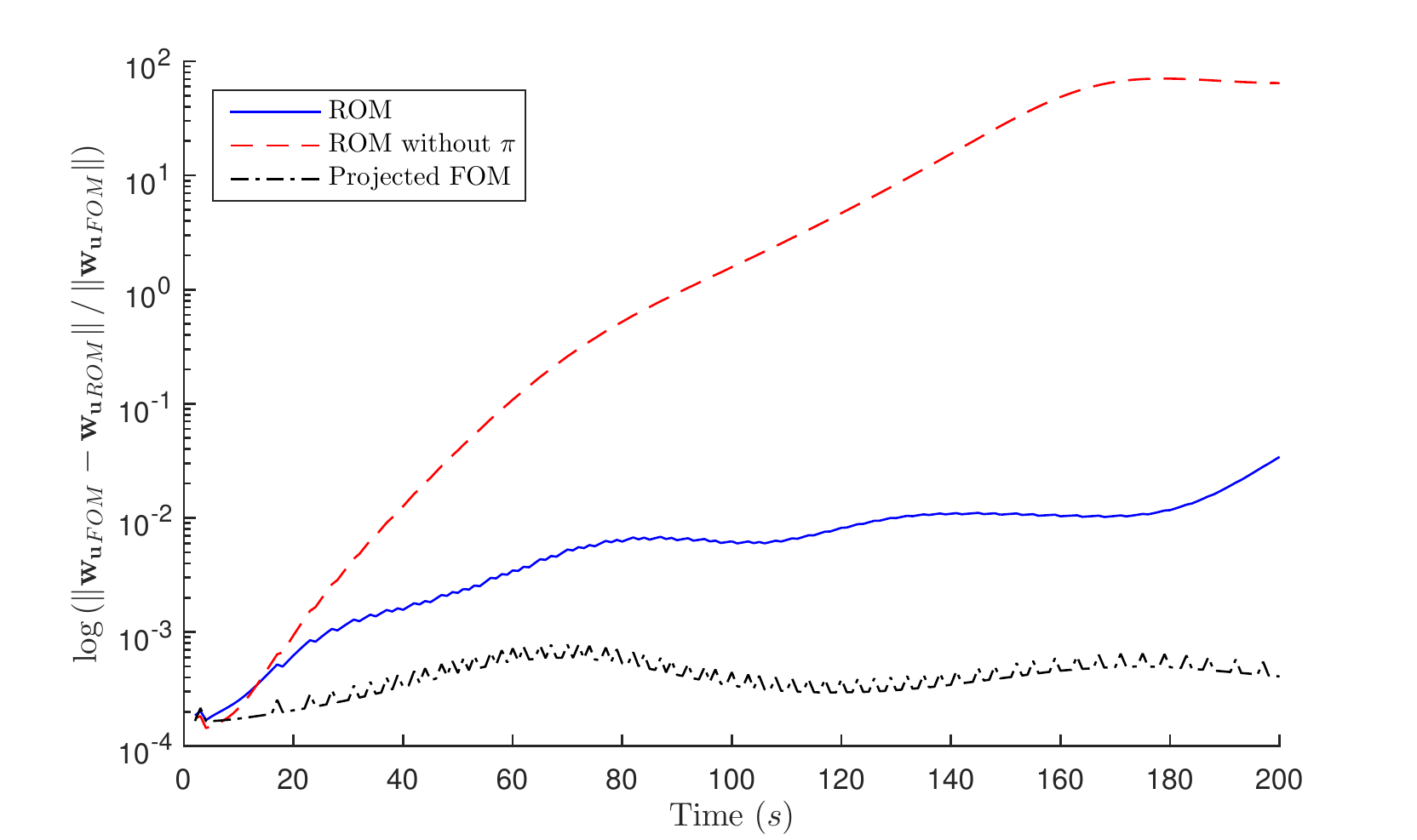}
	\caption{Manufactured test. Relative error of the linear momentum field for the projected FOM solution and the ROM solution with and without considering the pressure term in the momentum equation.}
	\label{fig:wu_t4n_250_sp}
\end{figure}

\begin{figure}[H]
	\centering
	\includegraphics[width=0.8\linewidth]{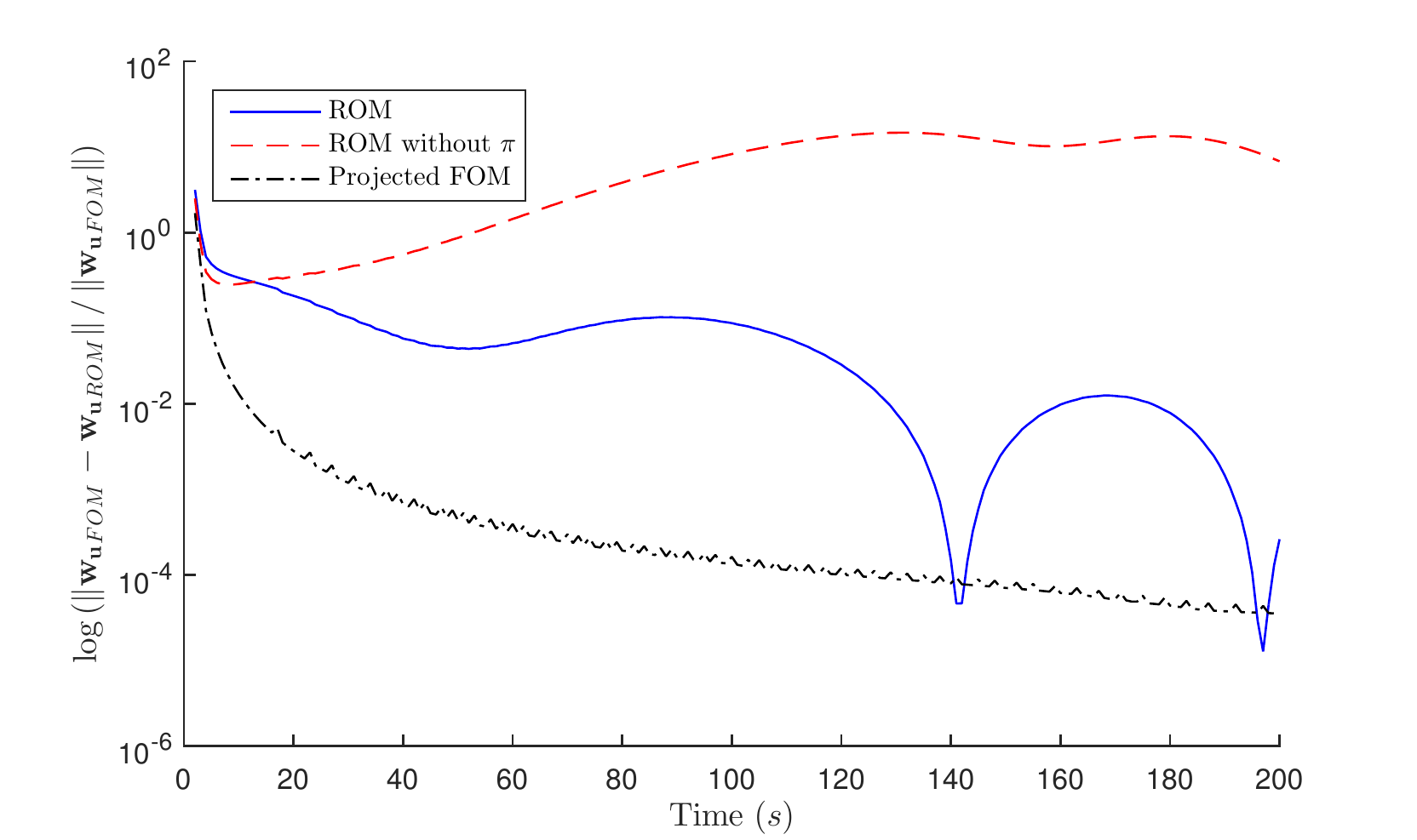}
	\caption{Manufactured test. Relative error of the pressure field for the projected FOM solution and the ROM solution with and without considering the pressure term in the momentum equation.}
	\label{fig:pi_t4n_250_sp}
\end{figure}

Let us remark that a key point on this test case is the consideration of pressure at the reduced order model. Neglecting the pressure gradient term in the momentum equation entails huge errors on the ROM solution. Figures \ref{fig:wu_t4n_250_sp} and \ref{fig:pi_t4n_250_sp} report the $L^{2}$ errors in logarithmic scale resulting from applying the ROM without the pressure term in equation \eqref{eq:dyn_sys_w}. Nevertheless, for some tests the magnitude of the pressure gradient term is considerably smaller than the one of the remaining terms in the equation and so, its computation could be avoided.

\subsection{Lid-driven cavity test with transport} \hfill\\

The second test is a modified version of the classical lid-driven cavity test in 3D by adding the resolution of a transport equation. The computational domain corresponds to the unit cube and the mesh is the one already employed at the first test.
The boundary of the domain is divided into two regions. At the top of the cavity we define an horizontal displacement  $\mathbf{u}=(1,0,0)^{T}$, whereas on the remaining boundaries we consider homogeneous Dirichlet boundary conditions for the velocity field. Moreover, Dirichlet homogeneous boundary conditions are defined for the species unknown on the whole boundary.
Regarding the initial conditions we assume zero velocity and 
we set the species to be given by
\begin{equation}
w_{\mathrm{y}}\left(\mathbf{x},0\right)=\left\lbrace \begin{array}{lr}
10 & \mathrm{if} \, (x_{1}-0.5)^{2}+(x_{2}-0.5)^{2}+(x_3-0.5)^{2}\leq 10^{-2},\\[4pt]
0 & \mathrm{if} \, (x_{1}-0.5)^{2}+(x_{2}-0.5)^{2}+(x_3-0.5)^{2}> 10^{-2}, 
\end{array}\right.
\end{equation}
that is, we define a ball of radio $10^{-2}$ in which the variable is set to $10$ and assume zero value on the remaining part of the domain. The density of the fluid is assumed to be constant and the viscosity is set to $\mu=10^{-2}$ and $\mathcal{D}=10^{-2}$.

The full order simulation is run up to time $t_{\mathrm{end}}=5$ with a fixed $CFL=1$. The snapshots are saved every $0.01$s. In the ROM different number of modes are considered attending to the fixed bounds $\kappa_{\mathbf{w}_{\mathbf{u}}}=99.99\%$,
$\kappa_{\pi}=99\%$ and
$\kappa_{w_{\mathrm{y}}}=99.99\%$ for the linear momentum, the pressure and the species respectively.
The cumulative eigenvalues are depicted in Table \ref{tab:ldc_cumulativeeigenvalues}.
A lifting function has been defined in order to homogenize the snapshots related to the linear momentum. Therefore, the number of elements of the basis of $X^{\mathrm{POD}}_{\mathbf{w_u}}$ is the number of modes determined by $\kappa_{\mathbf{w}_{\mathbf{u}}}$ plus one.
Accordingly, the dimensions of the basis are $N=8$, $N_{\pi}=2$ and $N_{w_{\mathrm{y}}}=9$. 
\begin{table}[h]
	\renewcommand{\arraystretch}{1.3}
	\begin{center}
		\begin{tabular}{|c||c|c|c|}
			\hline 
			Number of modes & $\mathbf{w}_{\mathbf{u}}$& $\pi$ & $w_{\mathrm{y}}$ \\ \hline
			\hline 
			$1$& $0.8681075$& $0.9792351$  & $0.7806067$ \\  
			$2$& $0.9742675$& $\mathbf{0.9980598}$  & $0.9456740$ \\ 
			$3$& $0.9944674$& $0.9998121$  & $0.9885488$ \\  
			$4$& $0.9988293$& $0.9999566$  & $0.9970488$ \\  
			$5$& $0.9996600$& $0.9999885$  & $0.9992604$ \\  
			$6$& $0.9998813$& $0.9999967$  & $0.9998145$ \\  
			$7$& $\mathbf{0.9999643}$& $0.9999983$  & $0.9999489$ \\  
			$8$& $0.9999854$& $0.9999990$  & $0.9999872$ \\
			$9$& $0.9999919$& $0.9999993$  & $\mathbf{0.9999968}$ \\ \hline 
		\end{tabular} 
	\end{center}
	
	\vspace*{0.25cm}
	\caption{Lid-driven cavity test. The second, third and fourth columns present the cumulative eigenvalues for the linear momentum, the pressure and the species, respectively. The values in which  the fixed bounds $\kappa_{\mathbf{w}_{\mathbf{u}}}\!=\!99.99\%$,
		$\kappa_{\pi}\!=\!99\%$ and
		$\kappa_{w_{\mathrm{y}}}\!=\!99.99\%$ are reached appear in bold font.}
	\label{tab:ldc_cumulativeeigenvalues}
\end{table}

\begin{figure}
	\begin{center}
		{\small \hspace*{1cm}	$t=0.1s$\hfill $t=0.5s$ \hfill $t=2.5s$ \hfill $t=5s$ \hspace*{1.2cm}}
	\end{center}	
	\includegraphics[width=0.24\linewidth]{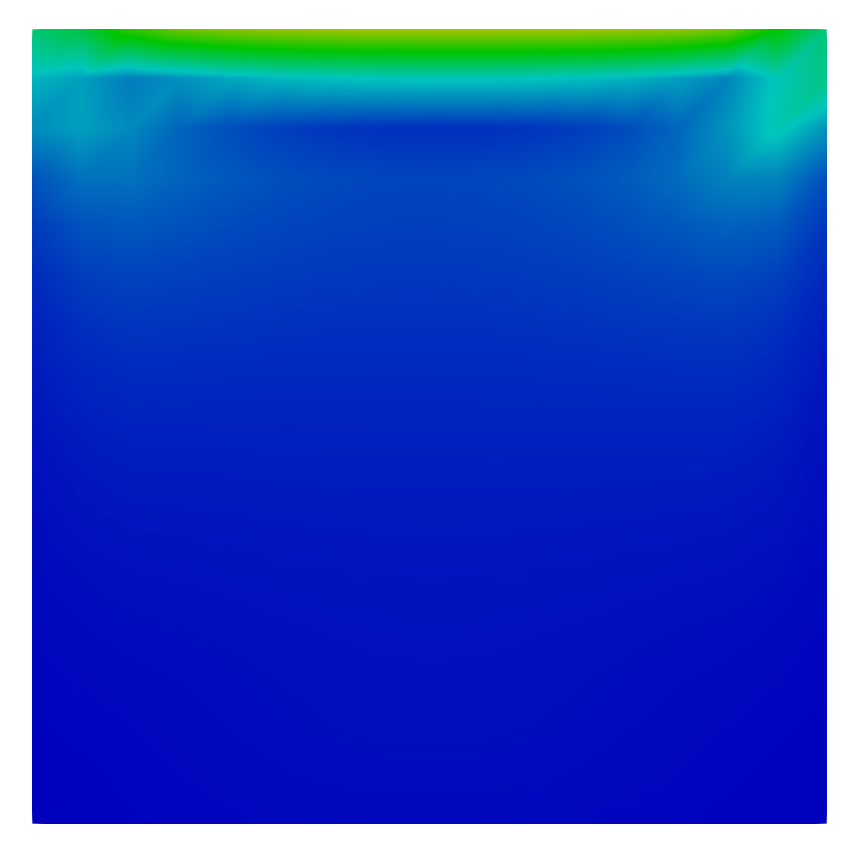}
	\hfill
	\includegraphics[width=0.24\linewidth]{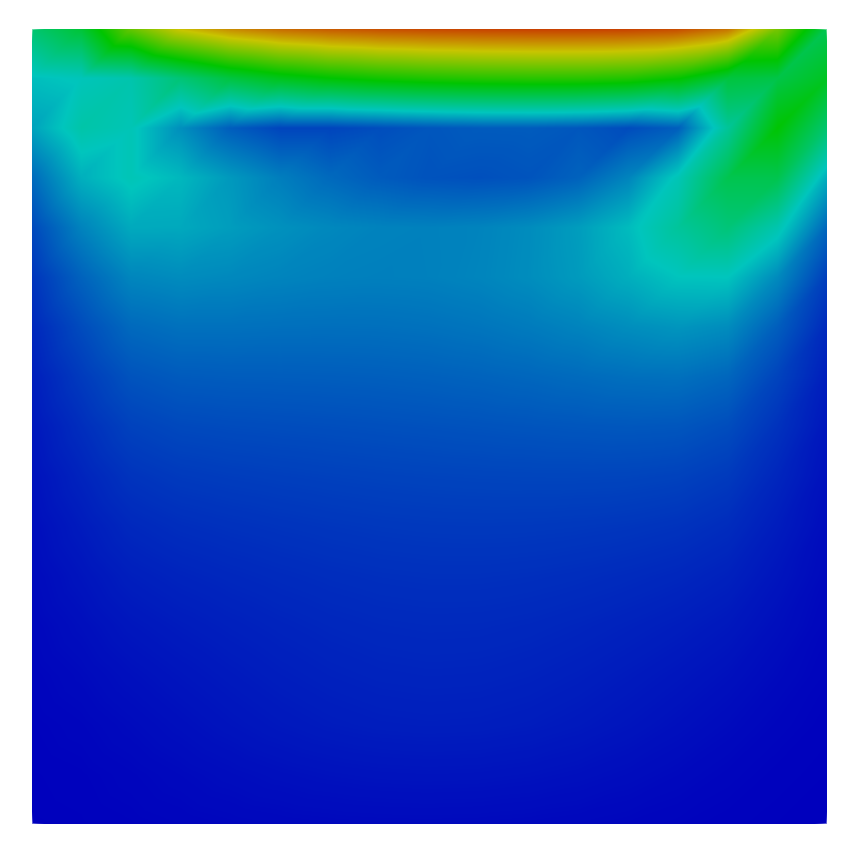}
	\hfill
	\includegraphics[width=0.24\linewidth]{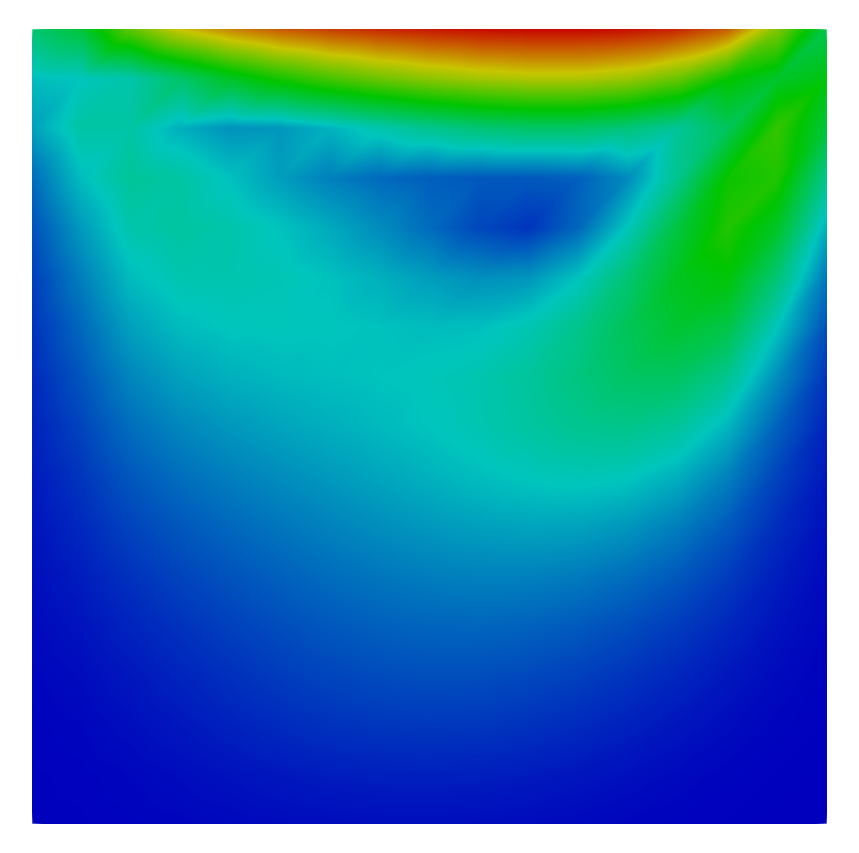}
	\hfill
	\includegraphics[width=0.24\linewidth]{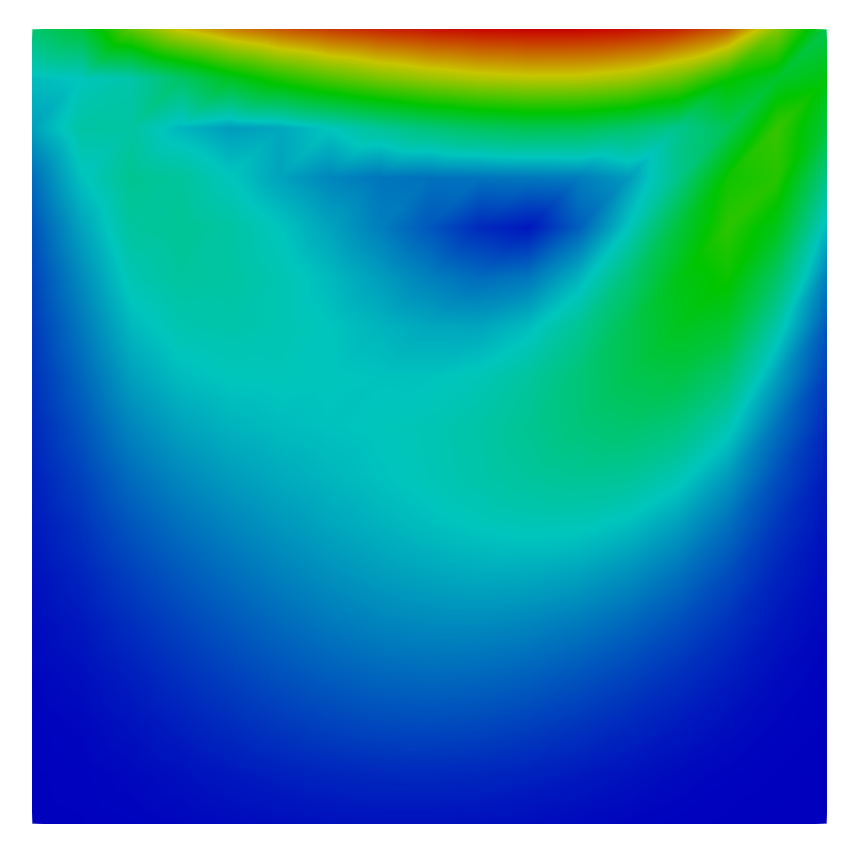}
	
	\includegraphics[width=0.24\linewidth]{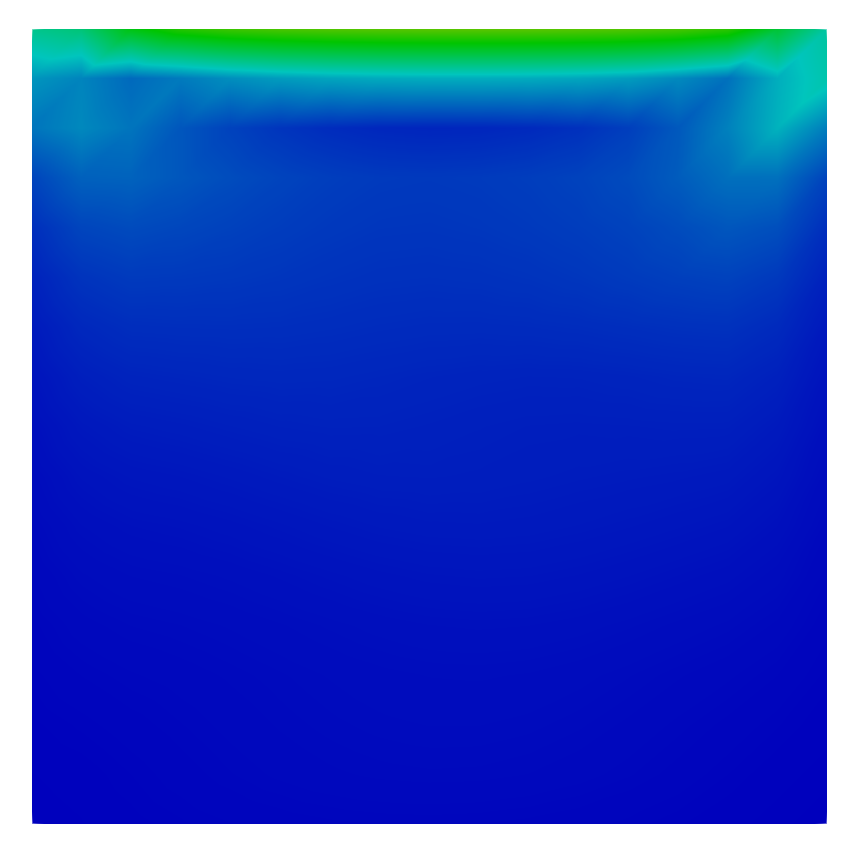}
	\hfill
	\includegraphics[width=0.24\linewidth]{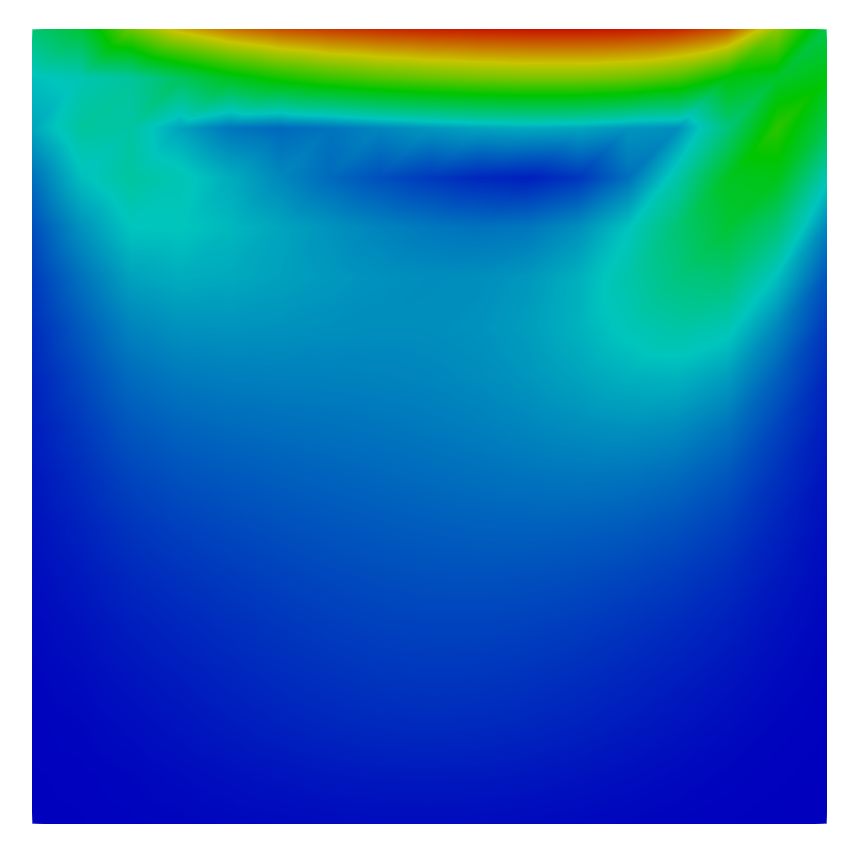}
	\hfill
	\includegraphics[width=0.24\linewidth]{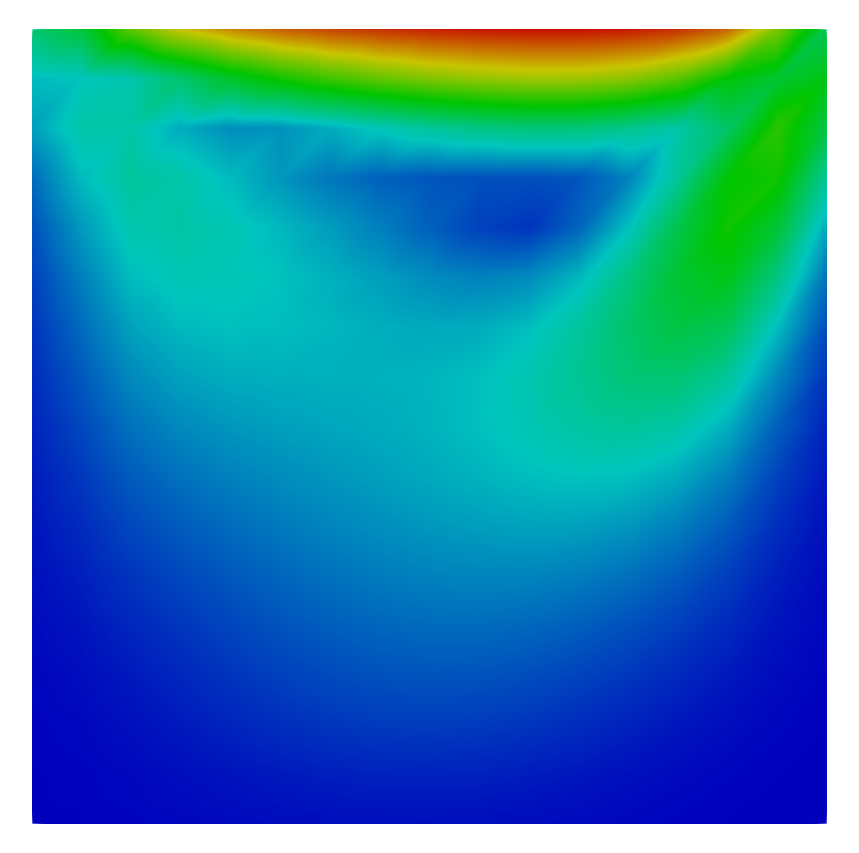}
	\hfill
	\includegraphics[width=0.24\linewidth]{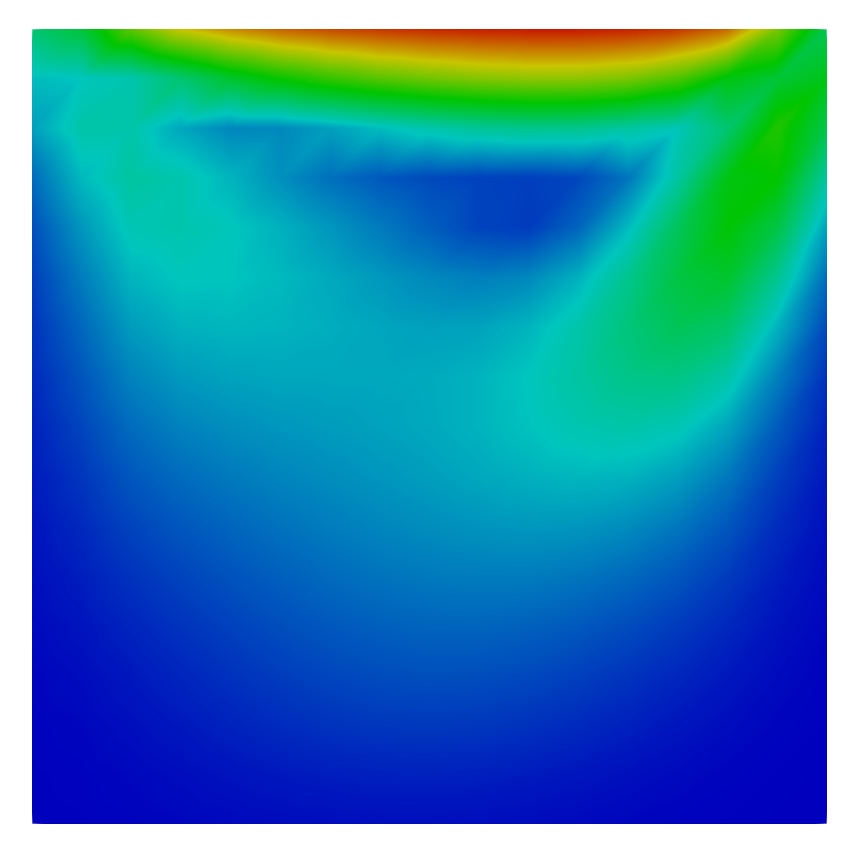}
	
	\vspace{0.2cm}
	\begin{center}
		\includegraphics[width=0.4\linewidth]{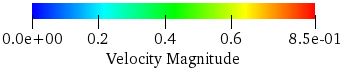}
	\end{center}
	
	\caption{Lid-driven cavity test. Comparison of the velocity Magnitude at plane $x=0.5$.
		Top: FOM solution. Bottom: ROM approximation.}
	\label{fig:ld16_5new10_velocitymagnitude}
\end{figure}

\begin{figure}
	\begin{minipage}{0.045\linewidth}		
	\end{minipage}\hfill
	\begin{minipage}{0.94\linewidth}
		\begin{center}
			\includegraphics[width=0.5\linewidth]{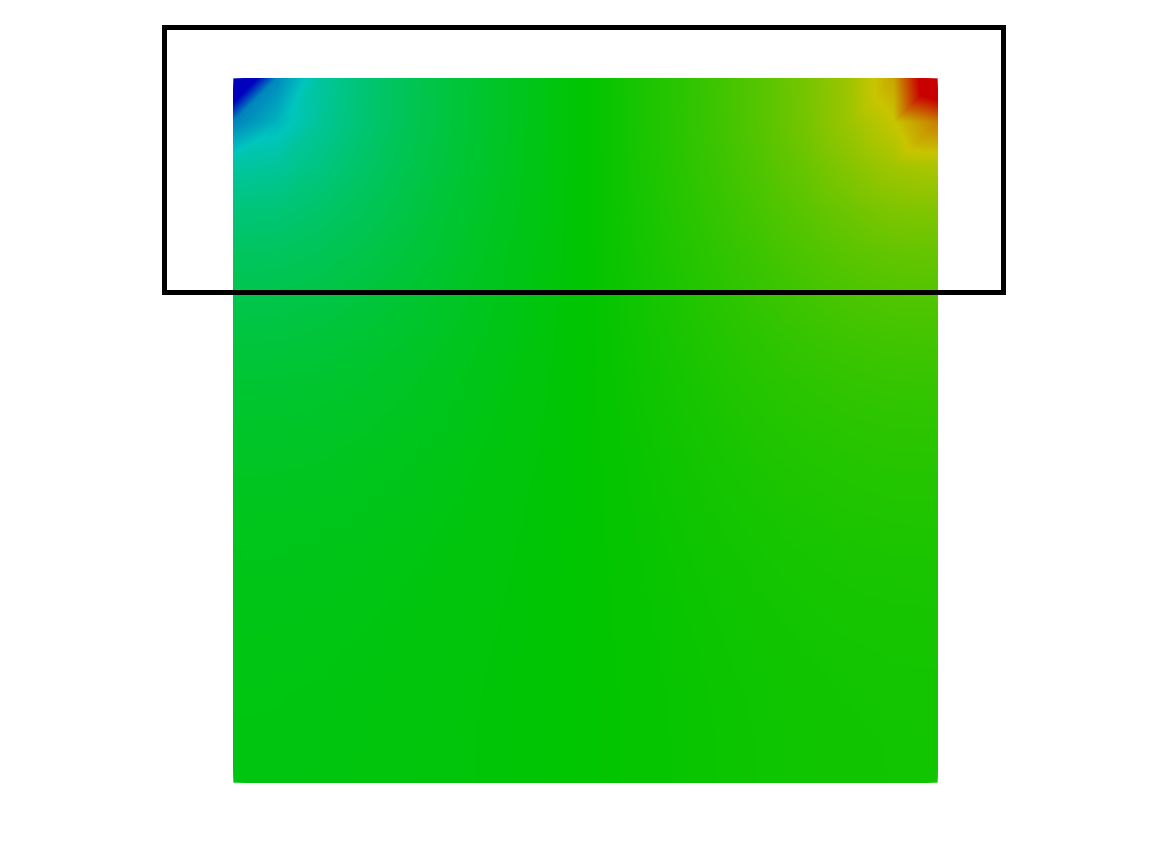}
	\end{center}\end{minipage}
	
	\begin{minipage}{0.1\linewidth}
		$t=0.1s$
	\end{minipage}\hspace*{-0.05\linewidth}
	\begin{minipage}{0.9\linewidth}
		\hfill\includegraphics[width=0.49\linewidth]{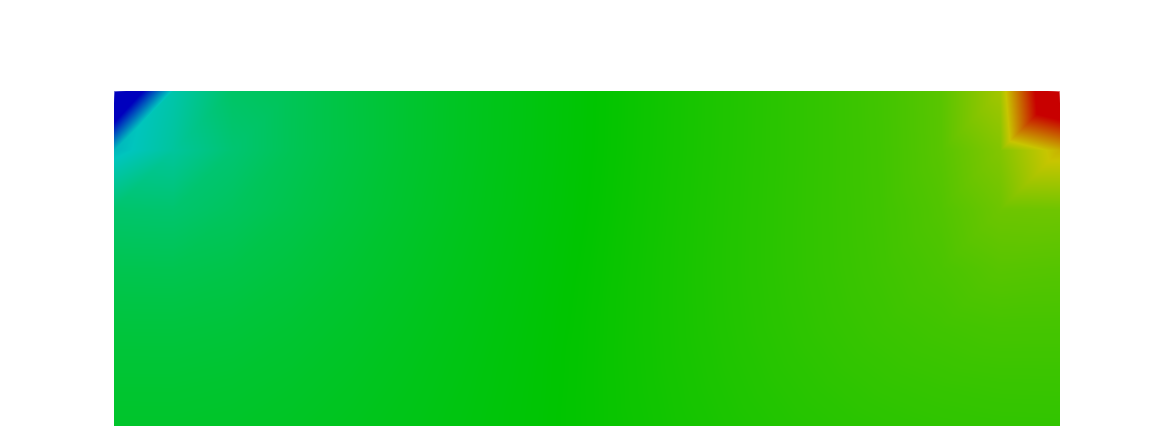}
		\hspace*{-0.05\linewidth}
		\includegraphics[width=0.49\linewidth]{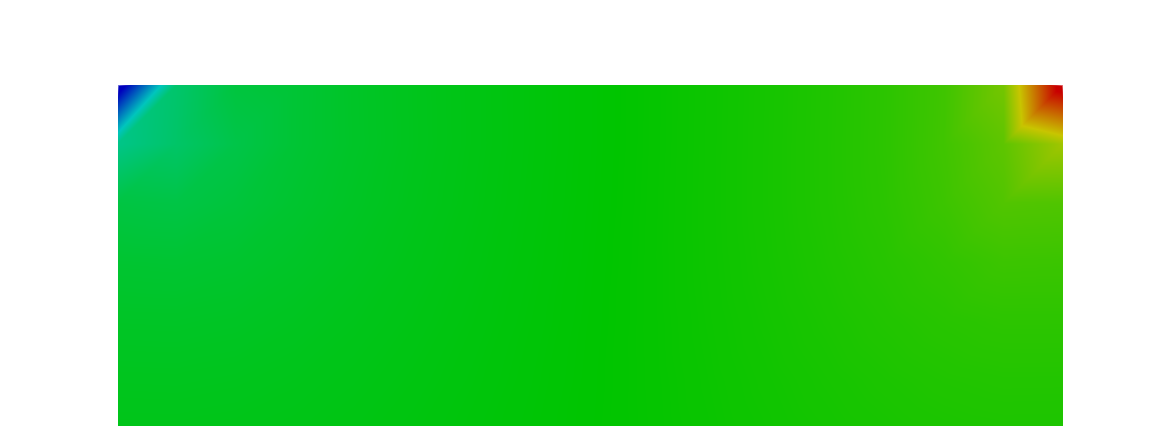}
	\end{minipage}
	
	\begin{minipage}{0.1\linewidth}
		$t=0.5s$
	\end{minipage}\hspace*{-0.05\linewidth}
	\begin{minipage}{0.9\linewidth}
		\hfill\includegraphics[width=0.49\linewidth]{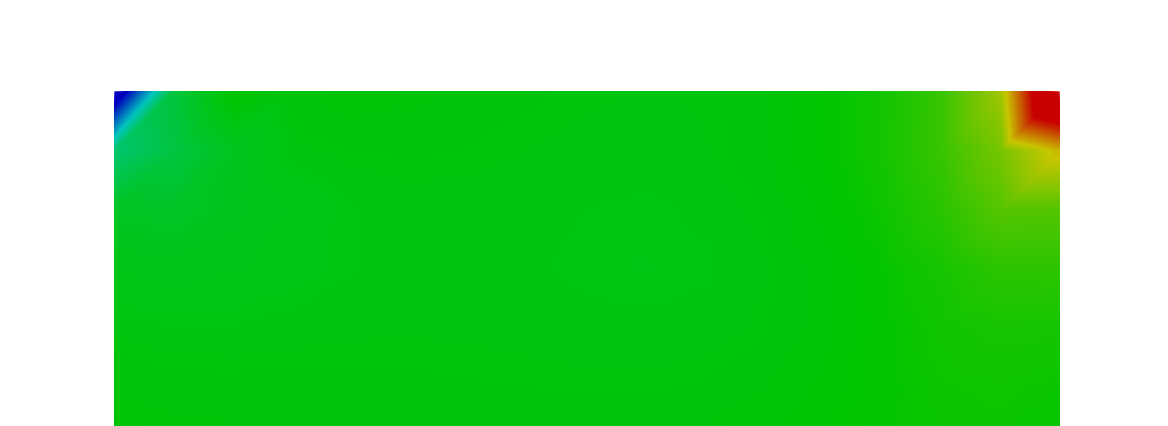}
		\hspace*{-0.05\linewidth}
		\includegraphics[width=0.49\linewidth]{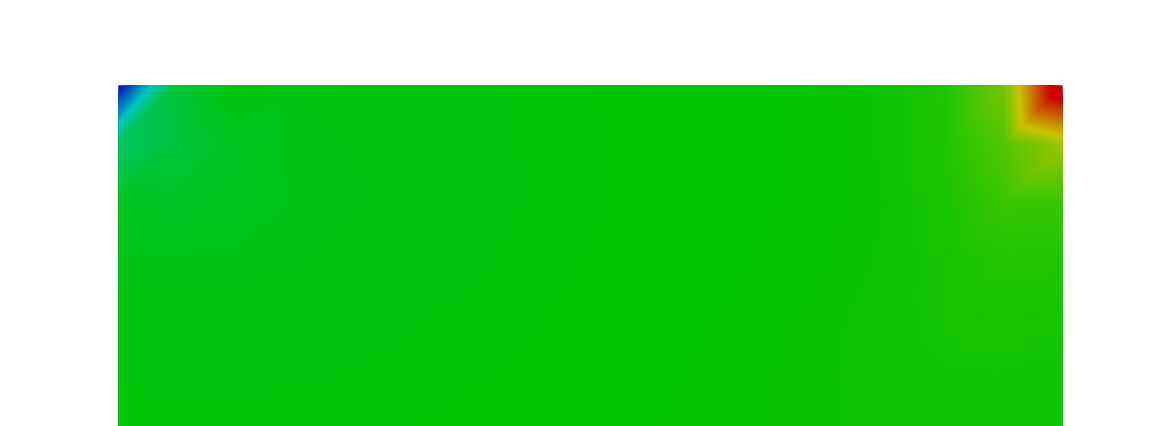}
	\end{minipage}
	
	\begin{minipage}{0.1\linewidth}
		$t=2.5s$
	\end{minipage}\hspace*{-0.05\linewidth}
	\begin{minipage}{0.9\linewidth} \hfill\includegraphics[width=0.49\linewidth]{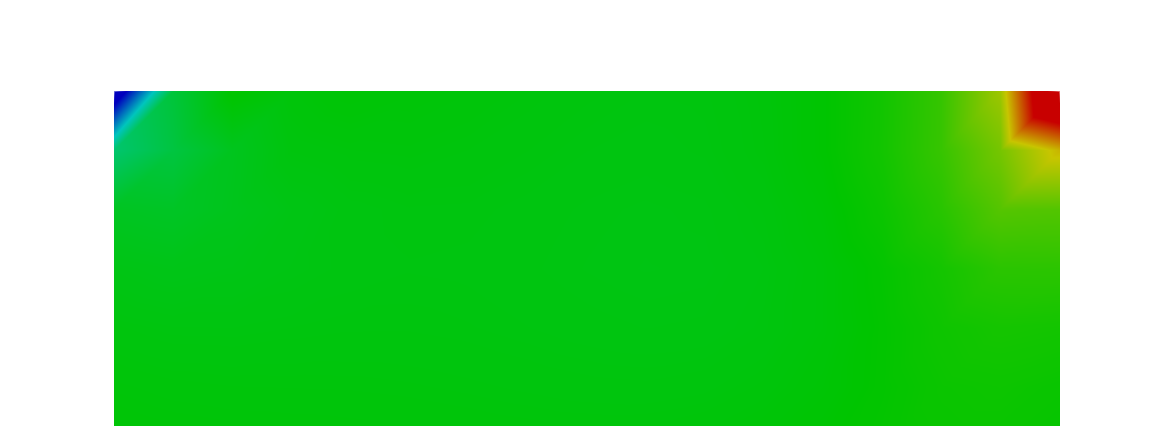}
		\hspace*{-0.05\linewidth}
		\includegraphics[width=0.49\linewidth]{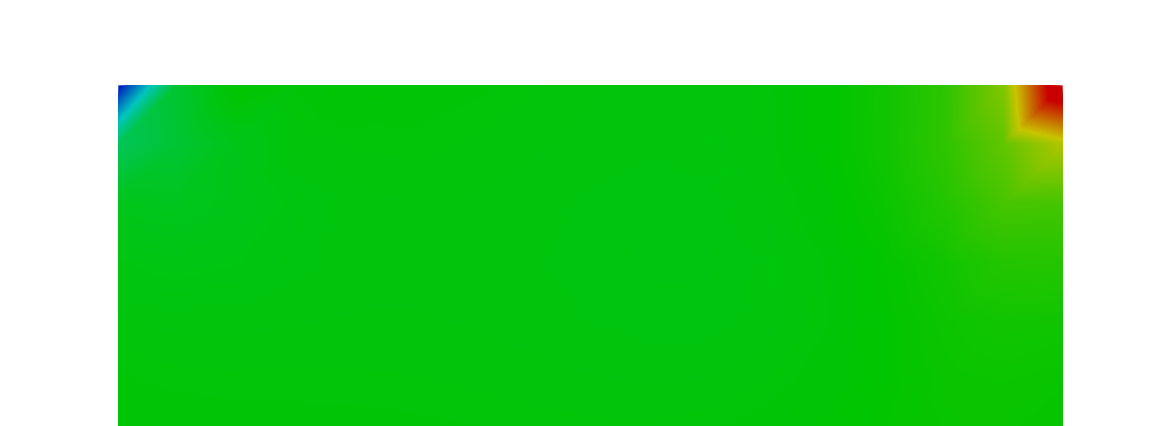}
	\end{minipage}
	
	\begin{minipage}{0.1\linewidth}
		$t=5s$
	\end{minipage}\hspace*{-0.05\linewidth}
	\begin{minipage}{0.9\linewidth} \hfill\includegraphics[width=0.49\linewidth]{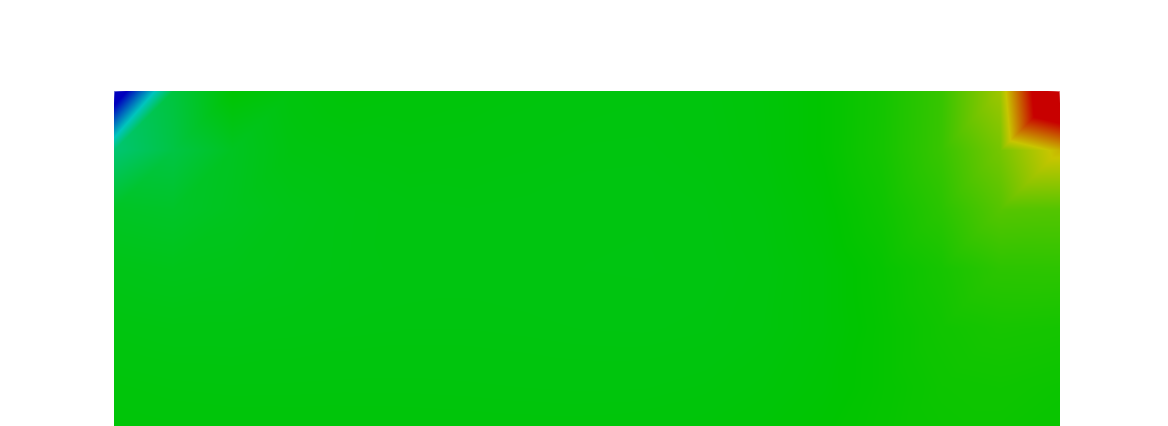}
		\hspace*{-0.05\linewidth}
		\includegraphics[width=0.49\linewidth]{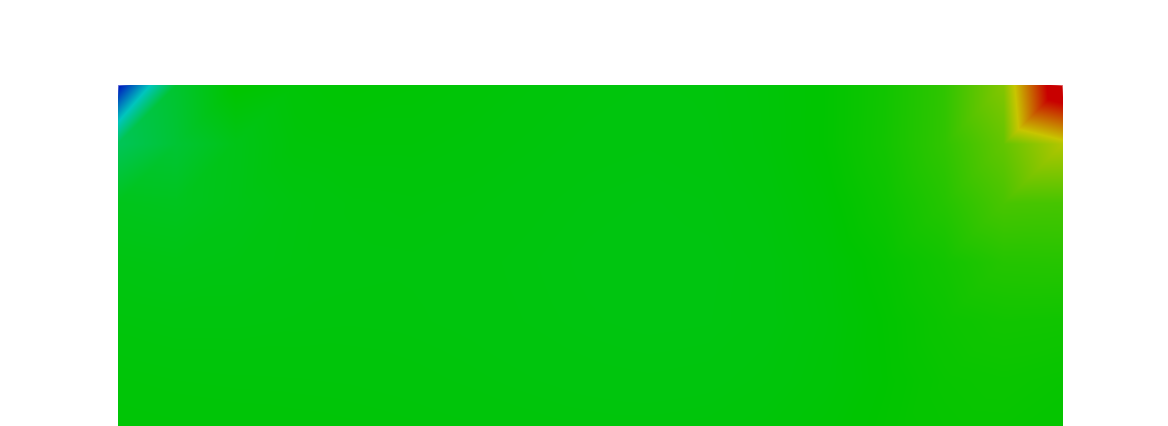}
	\end{minipage}
	
	\vspace{0.5cm}
	\begin{minipage}{0.05\linewidth}		
	\end{minipage}\hfill
	\begin{minipage}{0.94\linewidth}
		\begin{center}
			\includegraphics[width=0.4\linewidth]{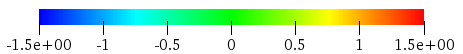}
		\end{center} 
	\end{minipage}

	\caption{Lid-driven cavity test. Comparison of the Pressure at plane $x=0.5$.
		Left: FOM solution. Right: ROM approximation.}
	\label{fig:ld16_5new10_pressure}
\end{figure}

\begin{figure}
	\begin{center}
		{\small \hspace*{1cm}	$t=0.1s$
			\hfill $t=1s$ \hfill $t=2.5s$ \hfill $t=5s$ \hspace*{1.2cm}}
	\end{center}	
	\includegraphics[width=0.24\linewidth]{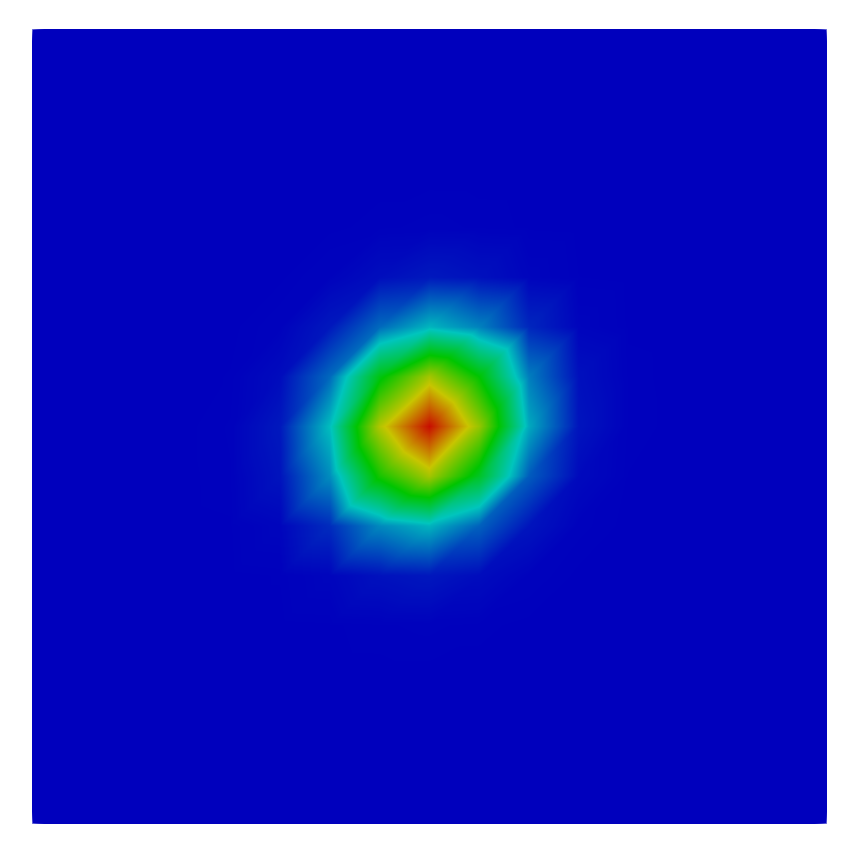}
	\hfill
	\includegraphics[width=0.24\linewidth]{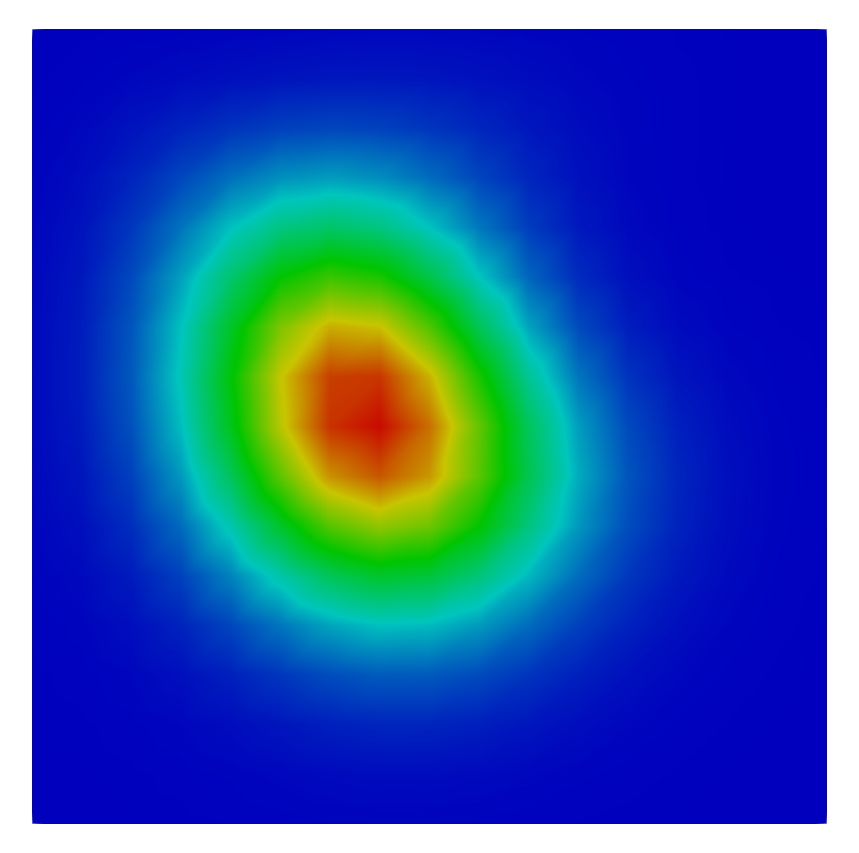}
	\hfill
	\includegraphics[width=0.24\linewidth]{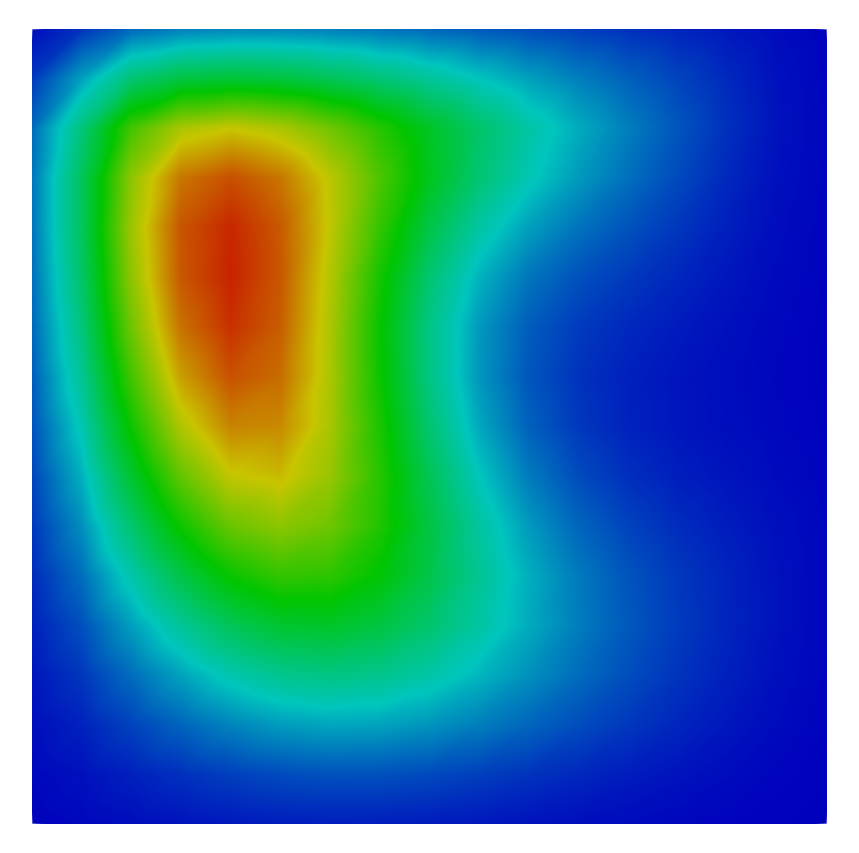}
	\hfill
	\includegraphics[width=0.24\linewidth]{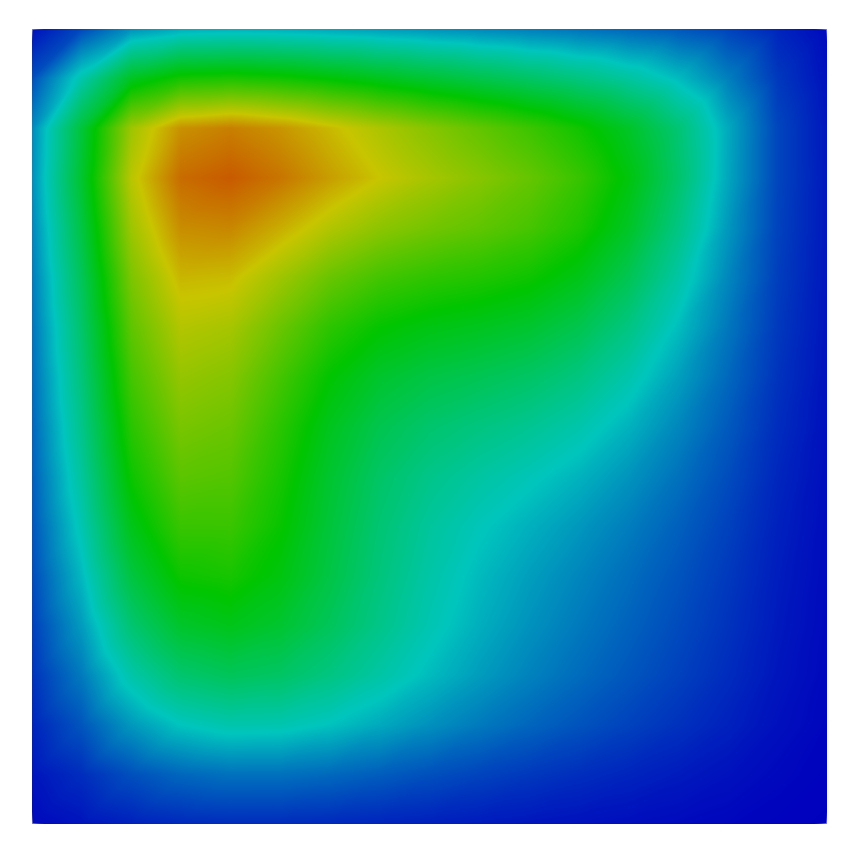}
	
	\includegraphics[width=0.24\linewidth]{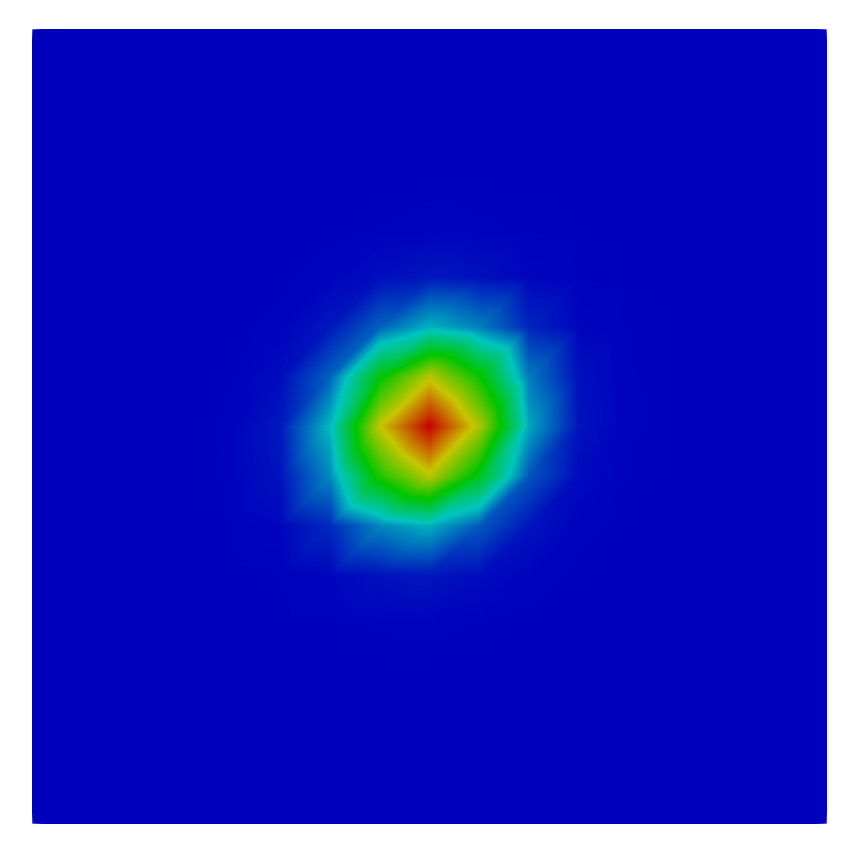}
	\hfill
	\includegraphics[width=0.24\linewidth]{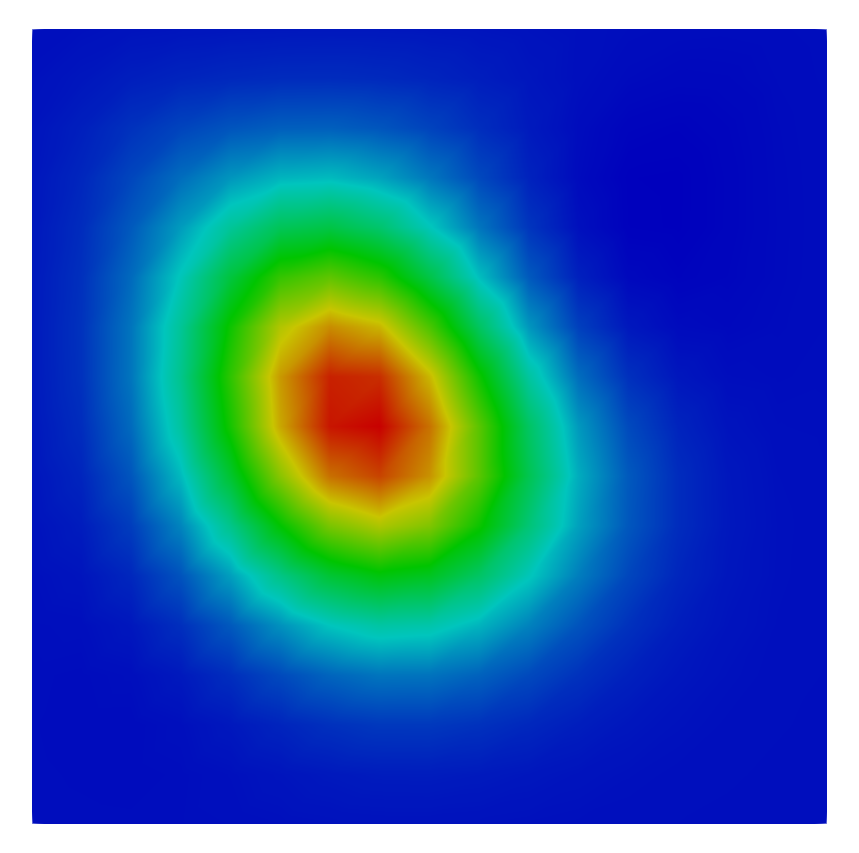}
	\hfill
	\includegraphics[width=0.24\linewidth]{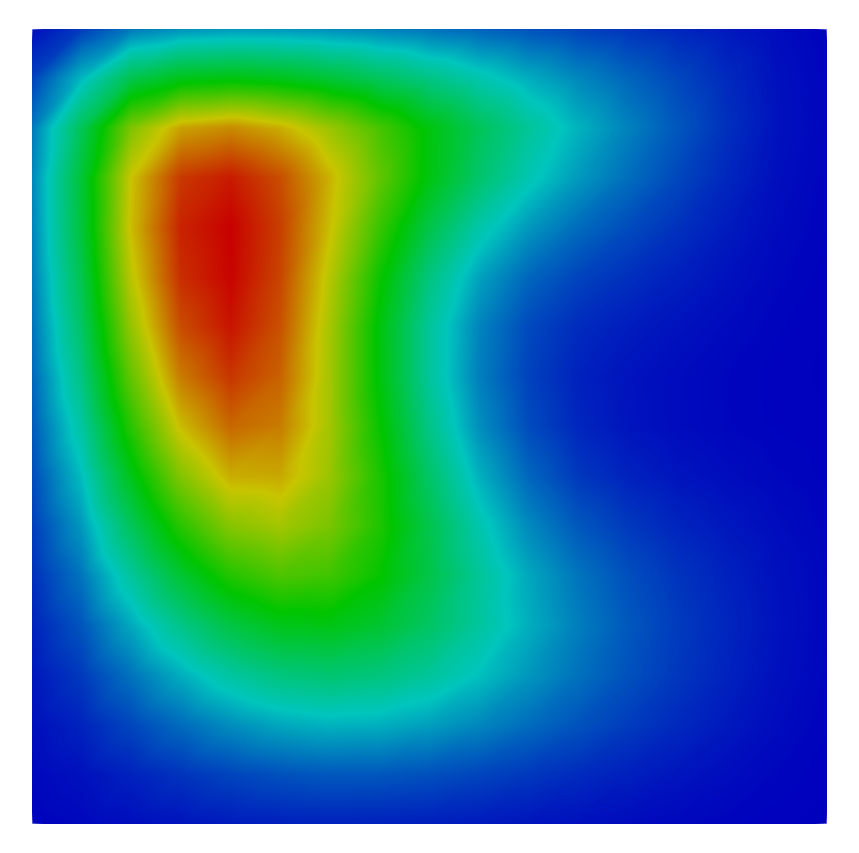}
	\hfill
	\includegraphics[width=0.24\linewidth]{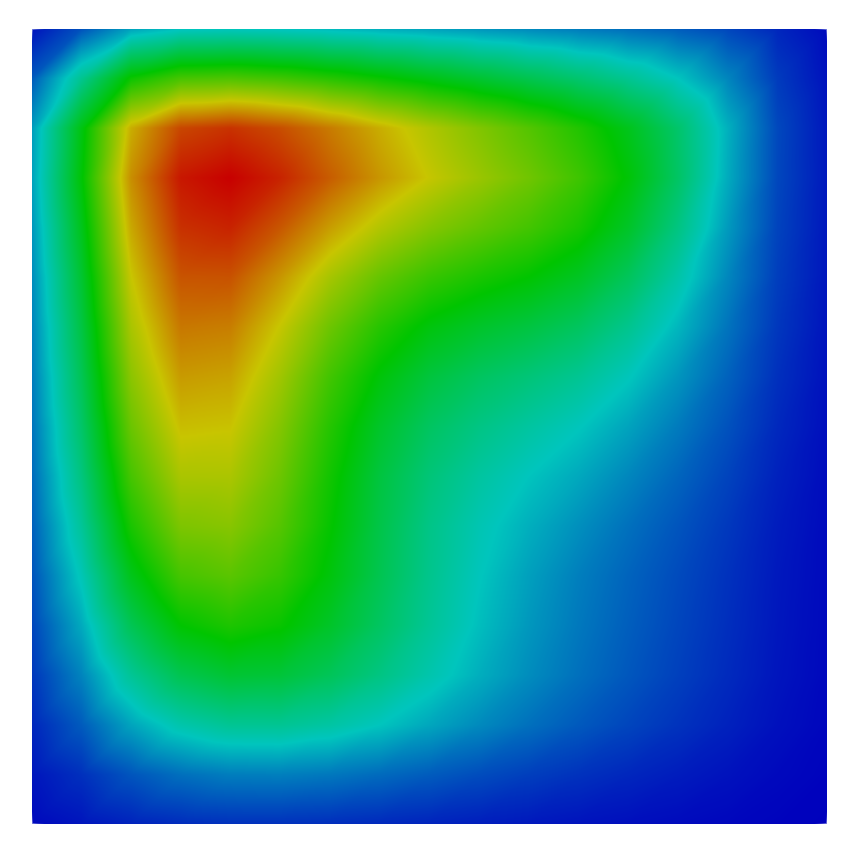}
	
	\vspace{0.2cm}
	\includegraphics[width=0.24\linewidth]{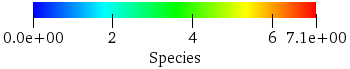}
	\hfill
	\includegraphics[width=0.24\linewidth]{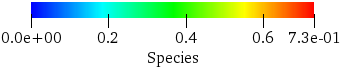}
	\hfill
	\includegraphics[width=0.24\linewidth]{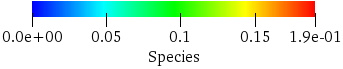}
	\hfill
	\includegraphics[width=0.24\linewidth]{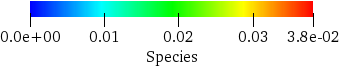}
	
	\caption{Lid-driven cavity test. Comparison of the Species at plane $x=0.5$.
		Top: FOM solution. Bottom: ROM approximation.}
	\label{fig:ld16_5new10_species}
\end{figure}

\begin{figure}
	\centering
	\includegraphics[width=0.8\linewidth]{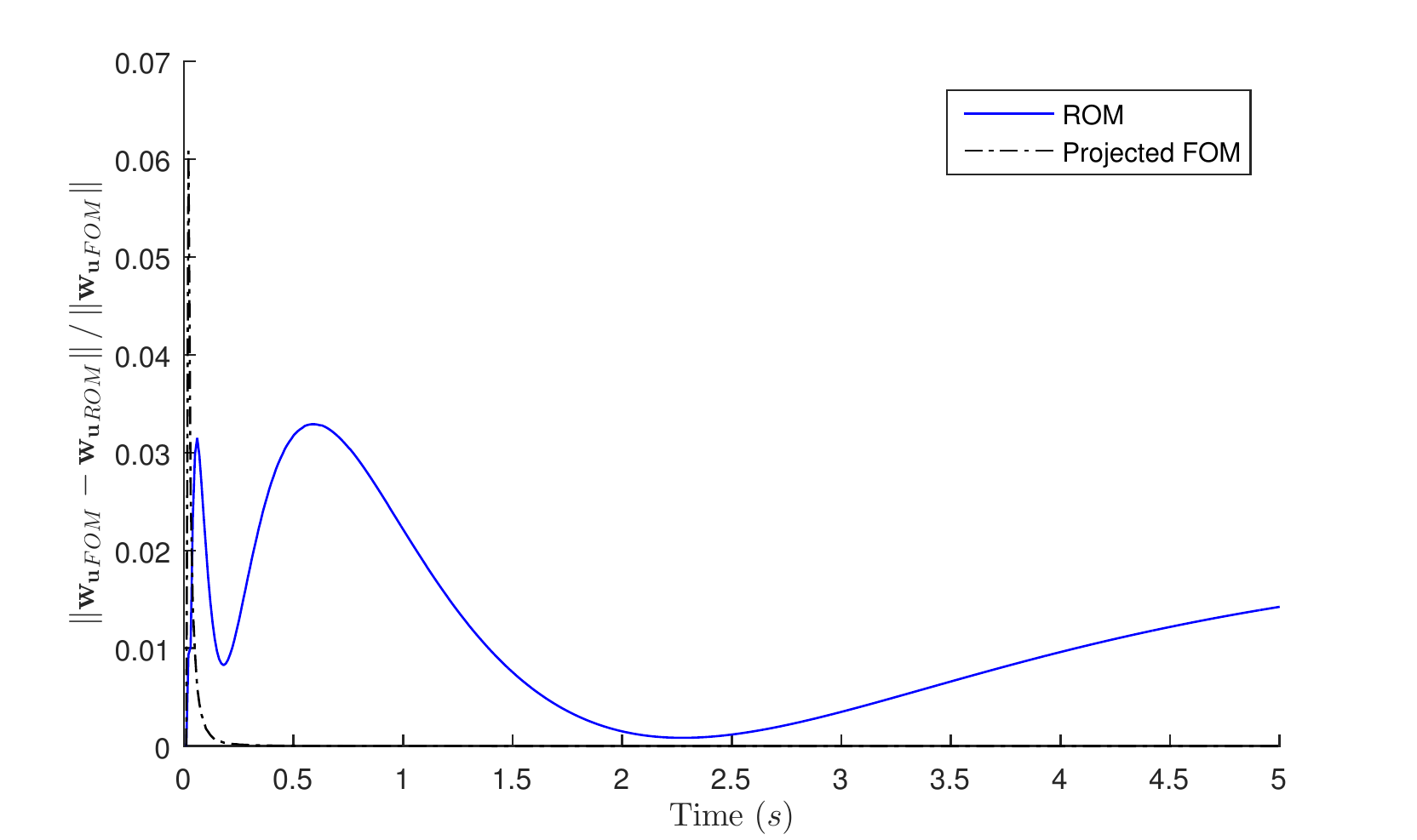}
	\caption{Lid-driven cavity test. Relative error of the linear momentum field for the projected FOM solution and the ROM solution. }
	\label{fig:wu10ball5newst}
\end{figure}

\begin{figure}
	\centering
	\includegraphics[width=0.8\linewidth]{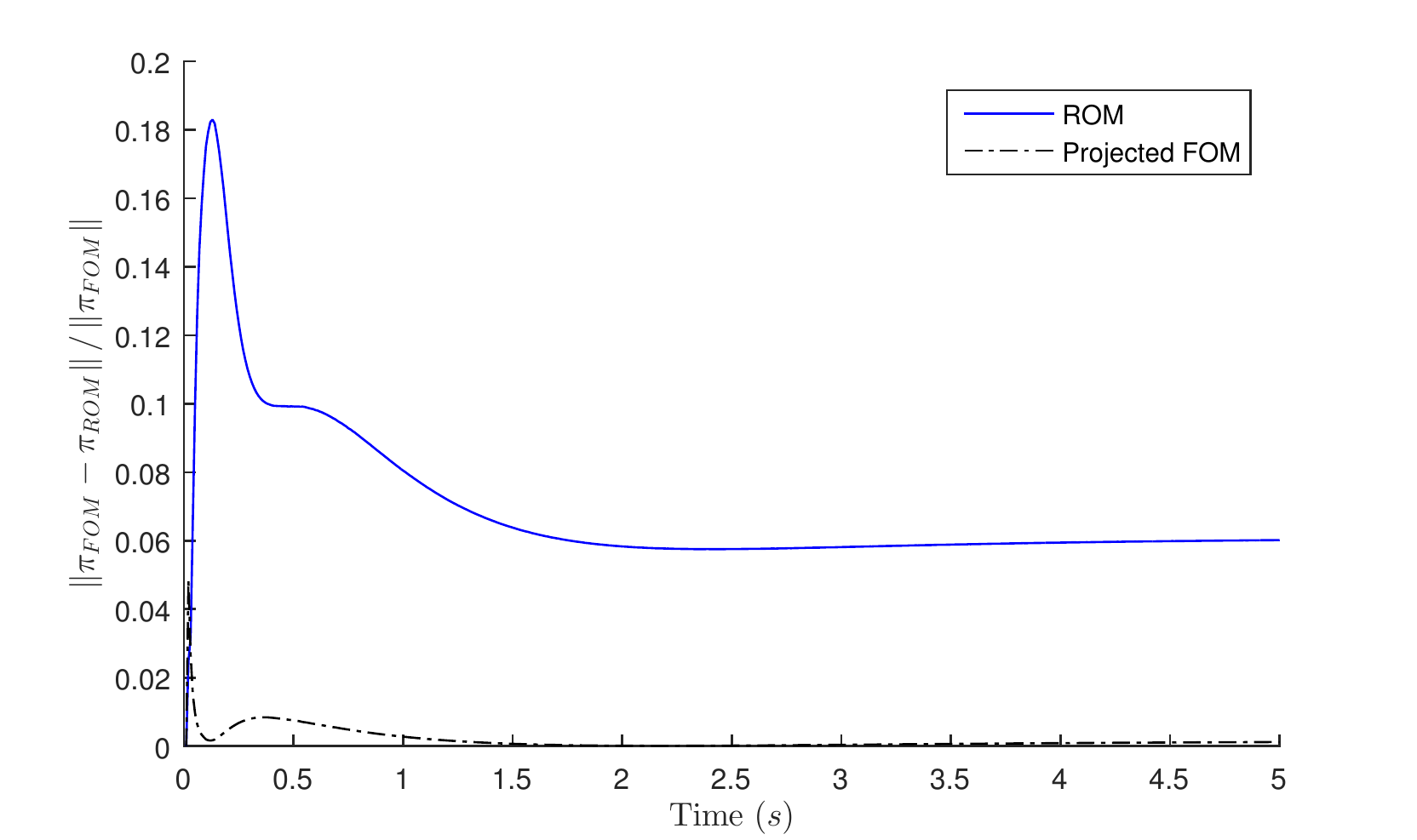}
	\caption{Lid-driven cavity test. Relative error of the pressure field for the projected FOM solution and the ROM solution. }
	\label{fig:pi10ball5newst}
\end{figure}

\begin{figure}
	\centering
	\includegraphics[width=0.8\linewidth]{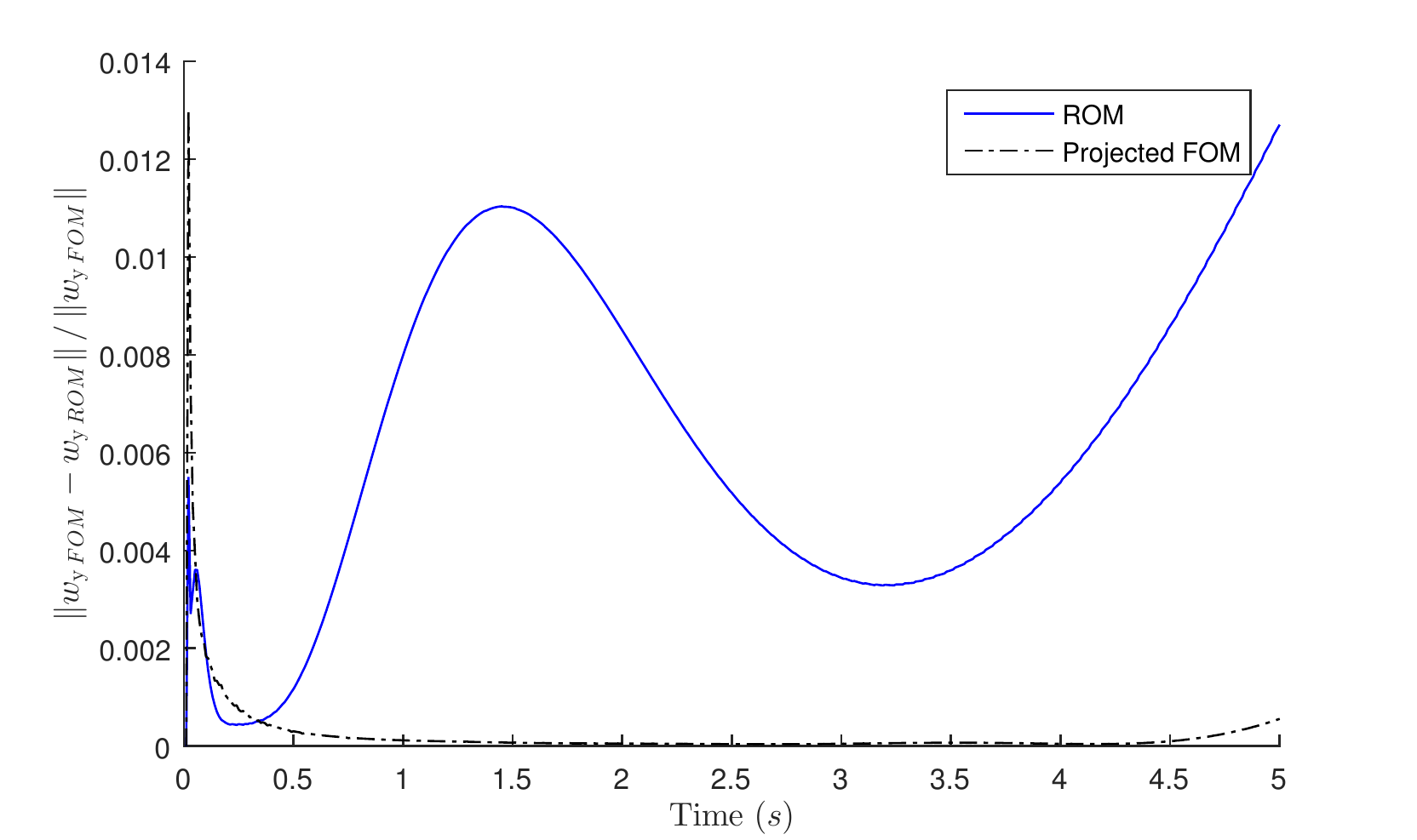}
	\caption{Lid-driven cavity test. Relative error of the species for the projected FOM solution and the ROM solution. }
	\label{fig:wy10ball5newst}
\end{figure}

The behaviour of the different fields as well as a comparison between FOM and ROM solutions is reported in Figures \ref{fig:ld16_5new10_velocitymagnitude}-\ref{fig:ld16_5new10_species}.
Figures \ref{fig:wu10ball5newst}-\ref{fig:wy10ball5newst} depict the obtained relative errors of the ROM solution and the projected solutions onto the selected POD basis. The results obtained show a good agreement with the FOM solution being able to capture the main effects of the flow. Nevertheless, several remarks should be taken into account.
Figure \ref{fig:wu10ball5newst} shows that the highest errors on the linear momentum field are found at initial times. This is due to the high variability of the velocity field at this period. Within the simulation, we have considered equally time-spaced snapshots. To enhance the performance of ROM, the concentration of snapshots should take into account the non-linear behaviour of the system. A similar issue can be observed in Figure \ref{fig:pi10ball5newst} for the pressure field.
On the other hand, the almost-linear behaviour of the pressure field has allowed us to consider a small number of basis elements to accurately reproduce the FOM solution. 
The fast decay of the species variable justifies the decrease on the accuracy graphically observed in the last column of Figure \ref{fig:ld16_5new10_species}. The magnitude of the species variable is slightly underestimated, whereas the features of it are properly captured. In this particular case enlarging the time interval, and therefore, the training set, for this variable would provide better results at $t=5$s.


\vspace{0.2cm}
% % % % % % % % % % % % % % % % % % % % % % % % % % % % % %
%                     Conclusions                         %
% % % % % % % % % % % % % % % % % % % % % % % % % % % % % %
\section{Conclusions and future developments}\label{sec:conclusions}

In this paper we have presented a novel POD-Galerkin reduced order method starting from an hybrid FV-FE full order solver for the incompressible Navier-Stokes equations in 3D. The method is based onto the definition of reduced basis functions defined on two staggered meshes (i.e the finite volume and the finite element one). 
The Galerkin projection of the momentum equation has been built to 
be consistent with the FV method used in the FOM taking advantage from the dual mesh structure. 
Special care has been taken in order to account for the pressure contribution in the ROM. The importance of including the pressure term onto the momentum equation has been proved by numerical results. Non-homogeneous boundary conditions have been overcome using a lifting function. 
The methodology has been verified on the unsteady Navier-Stokes equations giving promising results. 

Further, a  transport equation has been coupled with the incompressible Navier-Stokes equations and the corresponding modifications on the ROM have been presented. Besides, the good results obtained open the doors to the development of ROM for the resolution of turbulent flows by coupling the Navier-Stokes equations with the transport equations of a RANS or a LES turbulence model (\cite{LCLR16,enrique17,SBZ18,HSMR18,HSBR18}).

As future development, one of the main interest is into the analysis of the proposed methodology in a parametrized setting where, additionally to the time parameter, we will consider physical and geometrical parameters (\cite{ST18}).

 \vspace{0.2cm}
% % % % % % % % % % % % % % % % % % % % % % % % % % % % % %
%                  Acknowledgements                       %
% % % % % % % % % % % % % % % % % % % % % % % % % % % % % %
\section*{Acknowledgements}
We acknowledge the support provided by Spanish MECD under grant FPU13/00279, by Spanish MINECO under MTM2017-86459-R, by EU-COST MORNET TD13107 under STSM 40422, by the European Research Council Executive Agency with the Consolidator Grant project AROMA-CFD ``Advanced Reduced Order Methods with Applications in Computational Fluid Dynamics'' - GA 681447, H2020-ERC CoG 2015 AROMA-CFD (P.I. Gianluigi  Rozza) and by the INdAM-GNCS projects.

% % % % % % % % % % % % % % % % % % % % % % % % % % % % % %
%                    Bibliography                         %
% % % % % % % % % % % % % % % % % % % % % % % % % % % % % %
\bibliographystyle{amsplain}

\bibliography{./mibiblio}

\providecommand{\bysame}{\leavevmode\hbox to3em{\hrulefill}\thinspace}
\providecommand{\MR}{\relax\ifhmode\unskip\space\fi MR }
\providecommand{\MRhref}[2]{%
  \href{http://www.ams.org/mathscinet-getitem?mr=#1}{#2}
}
\providecommand{\href}[2]{#2}
\begin{thebibliography}{10}

\bibitem{ANR09}
I.~Akhtar, A.~H. Nayfeh, and C.~J. Ribbens, \emph{On the stability and
  extension of reduced-order {G}alerkin models in incompressible flows}, Theor.
  Comput. Fluid Dyn. \textbf{23} (2009), no.~3, 213--237.

\bibitem{Ball15}
F.~Ballarin, \emph{Reduced-order models for patient-specific haemodynamics of
  coronary artery bypass grafts}, Ph.D. thesis, 2015.

\bibitem{BMQR15}
F.~Ballarin, A.~Manzoni, A.~Quarteroni, and G.~Rozza, \emph{Supremizer
  stabilization of {POD-Galerkin} approximation of parametrized steady
  incompressible {Navier-Stokes} equations}, Int. J. Numer. Methods Eng.
  \textbf{102} (2015), no.~5, 1136--1161.

\bibitem{BARRAULT2004667}
M.~Barrault, Y.~Maday, N.C. Nguyen, and A.~Patera, \emph{An `empirical
  interpolation' method: application to efficient reduced-basis discretization
  of partial differential equations}, Comptes Rendus Mathematique \textbf{339}
  (2004), no.~9, 667--672.

\bibitem{bennerParSys}
P.~Benner, M.~Ohlberger, A.~T. Patera, G.~Rozza, and K.~Urban, \emph{Model
  reduction of parametrized systems}, vol. 1st 2017, MS\&A series, no.~17,
  Springer, 2017.

\bibitem{BBI09516}
M.~Bergmann, C.-H. Bruneau, and A.~Iollo, \emph{{Enablers for robust POD
  models}}, Journal of Computational Physics \textbf{228} (2009), no.~2,
  516--538.

\bibitem{BBFVC18}
A.~Berm\'udez, S.~Busto, J.~L. Ferr\'in, and M.~E. V\'azquez-Cend\'on, \emph{A
  high order projection method for low {M}ach number flows}, ArXiv preprint
  arXiv:1802.10585 (2018).

\bibitem{BDDV98}
A.~Berm\'udez, A.~Dervieux, J.~A. Desideri, and M.~E. V\'azquez-Cend\'on,
  \emph{Upwind schemes for the two-dimensional shallow water equations with
  variable depth using unstructured meshes}, Comput. Methods Appl. Mech. Eng.
  \textbf{155} (1998), no.~1, 49--72.

\bibitem{BFSV14}
A.~Berm\'udez, J.~L. Ferr\'in, L.~Saavedra, and M.~E. V\'azquez-Cend\'on,
  \emph{A projection hybrid finite volume/element method for low-{M}ach number
  flows}, J. Comp. Phys. \textbf{271} (2014), 360--378.

\bibitem{Boffi2013}
D.~Boffi, F.~Brezzi, and M.~Fortin, \emph{Mixed finite element methods and
  applications}, Springer Berlin Heidelberg, 2013.

\bibitem{BFTVC17}
S.~Busto, J.~L. Ferr\'in, E.~F. Toro, and M.~E. V\'azquez-Cend\'on, \emph{A
  projection hybrid high order finite volume/finite element method for
  incompressible turbulent flows}, J. Comput. Phys. \textbf{353} (2018),
  169--192.

\bibitem{BTVC16}
S.~Busto, E.~F. Toro, and M.~E. V\'azquez-Cend\'on, \emph{Design and analisis
  of {ADER}--type schemes for model advection--diffusion--reaction equations},
  J. Comp. Phys. \textbf{327} (2016), 553--575.

\bibitem{CIJS14598}
A.~Caiazzo, T.~Iliescu, V.~John, and S.~Schyschlowa, \emph{A numerical
  investigation of velocity-pressure reduced order models for incompressible
  flows}, J. Comput. Phys. \textbf{259} (2014), 598 -- 616.

\bibitem{Carlberg2013623}
K.~Carlberg, C.~Farhat, J.~Cortial, and D.~Amsallem, \emph{The {GNAT} method
  for nonlinear model reduction: Effective implementation and application to
  computational fluid dynamics and turbulent flows}, J. Comput. Phys.
  \textbf{242} (2013), 623--647.

\bibitem{CVC10}
L.~Cea and M.~E. V{\'a}zquez-Cend{\'o}n, \emph{Unstructured finite volume
  discretization of two-dimensional depth-averaged shallow water equations with
  porosity}, Int. J. Numer. Methods Fluids \textbf{63} (2010), no.~8, 903--930.

\bibitem{CV12}
L.~Cea and M.~E. V\'azquez-Cend\'on, \emph{Analysis of a new {K}olgan-type
  scheme motivated by the shallow water equations}, Appl. Num. Math.
  \textbf{62} (2012), no.~4, 489--506.

\bibitem{ChinestaEnc2017}
F.~Chinesta, A.~Huerta, G.~Rozza, and K.~Willcox, \emph{Model order reduction},
  Encyclopedia of Computational Mechanics, Elsevier Editor, 2016 (2016).

\bibitem{DKKO91}
A.~E. Deane, I.~G. Kevrekidis, G.~E. Karniadakis, and S.~A. Orszag,
  \emph{Low-dimensional models for complex geometry flows: Application to
  grooved channels and circular cylinders}, Phys. Fluids \textbf{3} (1991),
  no.~10, 2337--2354.

\bibitem{DHO2012}
M.~Drohmann, B.~Haasdonk, and M.~Ohlberger, \emph{Reduced basis approximation
  for nonlinear parametrized evolution equations based on empirical operator
  interpolation}, SIAM J. Sci. Comput. \textbf{34} (2012), no.~2, A937--A969.

\bibitem{Everson1995}
R.~Everson and L.~Sirovich, \emph{Karhunen-{L}o{\`{e}}ve procedure for gappy
  data}, Journal of the Optical Society of America A \textbf{12} (1995), no.~8,
  1657.

\bibitem{Graham1999}
W.~R. Graham, J.~Peraire, and K.~Y. Tang, \emph{Optimal control of vortex
  shedding using low-order models. {P}art {II}-model-based control},
  International Journal for Numerical Methods in Engineering \textbf{44}
  (1999), no.~7, 973--990.

\bibitem{GS87}
P.~M. Gresho and R.~L. Sani, \emph{On pressure boundary conditions for the
  incompressible {N}avier-{S}tokes equations}, Int. J. Numer. Methods Fluids
  \textbf{7} (1987), no.~10, 1111--1145.

\bibitem{Guer06}
J.~L. Guermond, P.~Minev, and J.~Shen, \emph{An overview of projection methods
  for incompressible flows}, Comput. Methods Appl. Mech. Eng. \textbf{195}
  (2006), 6011--6045.

\bibitem{Gunzburger2007}
M.~D. Gunzburger, J.~S. Peterson, and J.~N. Shadid, \emph{Reduced-order
  modeling of time-dependent {PDEs} with multiple parameters in the boundary
  data}, Comput. Methods Appl. Mech. Eng. \textbf{196} (2007), no.~4-6,
  1030--1047.

\bibitem{HOR08}
B.~Haasdonk, M.~Ohlberger, and G.~Rozza, \emph{A reduced basis method for
  evolution schemes with parameter-dependent explicit operators}, ETNA,
  Electronic Transactions on Numerical Analysis \textbf{32} (2008), 145--161.

\bibitem{HRS16}
J.~S. Hesthaven, G.~Rozza, and B.~Stamm, \emph{Certified reduced basis methods
  for parametrized partial differential equations}, Springer, 2016.

\bibitem{HSBR18}
S.~Hijazi, S.~Ali, G.~Stabile, F.~Ballarin, and G.~Rozza, \emph{The effort of
  increasing {R}eynolds number in projection-based reduced order methods: from
  laminar to turbulent flows}, In press, FEF special volume 2017.

\bibitem{IL98a}
Angelo Iollo and St{\'e}phane Lanteri, \emph{{Approximation of compressible
  flows by a reduced order model}}, pp.~55--60, Springer Berlin Heidelberg,
  Berlin, Heidelberg, 1998.

\bibitem{IR98}
K.~Ito and S.S. Ravindran, \emph{{A Reduced-Order Method for Simulation and
  Control of Fluid Flows}}, Journal of Computational Physics \textbf{143}
  (1998), no.~2, 403 -- 425.

\bibitem{JL04}
H.~Johnston and J.-G. Liu, \emph{Accurate, stable and efficient
  {N}avier-{S}tokes solvers based on explicit treatment of the pressure term},
  J. Comput. Phys. \textbf{199} (2004), no.~1, 221--259.

\bibitem{Pet89}
Peterson. J.S., \emph{{The Reduced Basis Method for Incompressible Viscous Flow
  Calculations}}, SIAM Journal on Scientific and Statistical Computing
  \textbf{10} (1989), no.~4, 777--786.

\bibitem{KV02}
K.~Kunisch and S.~Volkwein, \emph{Galerkin proper orthogonal decomposition
  methods for a general equation in fluid dynamics}, SIAM J Numer Anal
  \textbf{40} (2002), no.~2, 492--515.

\bibitem{LLP10}
J.-G. Liu, J.~Liu, and R.~L. Pego, \emph{Stable and accurate pressure
  approximation for unsteady incompressible viscous flow}, J. Comput. Phys.
  \textbf{229} (2010), no.~9, 3428--3453.

\bibitem{LCLR16}
S.~Lorenzi, A.~Cammi, L.~Luzzi, and G.~Rozza, \emph{{POD-Galerkin} method for
  finite volume approximation of {Navier-Stokes} and {RANS} equations},
  Computer Methods in Applied Mechanics and Engineering \textbf{311} (2016),
  151 -- 179.

\bibitem{NE94}
B.~R. Noack and H.~Eckelmann, \emph{A low-dimensional {G}alerkin method for the
  three-dimensional flow around a circular cylinder}, Physics of Fluids
  \textbf{6} (1994), no.~1, 124--143.

\bibitem{NPP05}
Bernd~R. Noack, Paul Papas, and Peter~A. Monkewitz, \emph{The need for a
  pressure-term representation in empirical {G}alerkin models of incompressible
  shear flows}, Journal of Fluid Mechanics \textbf{523} (2005), 339--365.

\bibitem{OID86}
S.~A. Orszag, M.~Israeli, and M.~O. Deville, \emph{Boundary conditions for
  incompressible flows}, J. Sci. Comput. \textbf{1} (1986), no.~1, 75--111.

\bibitem{quarteroniRB2016}
A.~Quarteroni, A.~Manzoni, and F.~Negri, \emph{Reduced basis methods for
  partial differential equations}, Springer International Publishing, 2016.

\bibitem{enrique17}
T.~Chac\'on Rebollo, E.~Delgado \'Avila, M.~G\'omez M\'armol, F.~Ballarin, and
  G.~Rozza, \emph{On a certified {S}magorinsky reduced basis turbulence model},
  SIAM J. Numer. Anal. \textbf{55} (2017), no.~6, 3047–3067.

\bibitem{Rozza2007}
G.~Rozza and K.~Veroy, \emph{On the stability of the reduced basis method for
  {S}tokes equations in parametrized domains}, Comput. Methods Appl. Mech. Eng.
  \textbf{196} (2007), no.~7, 1244--1260.

\bibitem{HSMR18}
S.Hijazi, G.~Stabile, A.~Mola, and G.~Rozza, \emph{Data-driven {POD-G}alerkin
  reduced order model for turbulent flows}, In preparation (2018).

\bibitem{Sirisup2005}
S.~Sirisup and G.~E. Karniadakis, \emph{Stability and accuracy of periodic flow
  solutions obtained by a {POD}-penalty method}, Physica D: Nonlinear Phenomena
  \textbf{202} (2005), no.~3-4, 218--237.

\bibitem{SBZ18}
G.~Stabile, F.~Ballarin, G.~Zuccarino, and G.~Rozza, \emph{A reduced order
  variational multiscale approach for turbulent flows}, Submitted (2018).

\bibitem{SHMLR17}
G.~Stabile, S.~Hijazi, A.~Mola, S.~Lorenzi, and G.~Rozza,
  \emph{{POD}-{G}alerkin reduced order methods for {CFD} using finite volume
  discretisation: vortex shedding around a circular cylinder}, Communications
  in Applied and Industrial Mathematics \textbf{8} (2017).

\bibitem{SR17}
G.~Stabile and G.~Rozza, \emph{Finite volume {POD-G}alerkin stabilised reduced
  order methods for the parametrised incompressible {N}avier�{S}tokes
  equations}, Comput. Fluids \textbf{173} (2018), 273--284.

\bibitem{ST18}
\bysame, \emph{Geometrical parametrization for finite-volume based reduced
  order methods}, 2018.

\bibitem{TMN01}
E.~F. Toro, R.~C. Millington, and L.~A.~M. Nejad, \emph{Godunov methods},
  ch.~Towards very high order {G}odunov schemes, Springer, 2001.

\end{thebibliography}
\end{document}